\documentclass{amsart}

\usepackage{amscd}
\usepackage{verbatim}
\usepackage{amssymb}

\theoremstyle{plain}
\newtheorem{Thm}{Theorem}
\newtheorem{Cor}{Corollary}

\newtheorem{Prop}{Proposition}

\theoremstyle{definition}
\newtheorem{Def}{Definition}

\theoremstyle{remark}

\numberwithin{equation}{section}
\newcommand{\PS}[1]{{\mathbb P} ^{#1}}                  
\newcommand{\Op}[1]{{\mathcal O}_{{\mathbb P}^{#1}}} 
\newcommand{\sexs}[3]{0 \longrightarrow #1 \longrightarrow #2
\longrightarrow #3 \longrightarrow 0}
\newcommand{\res}[1]{\!\!\mid_{#1}}
\newcommand{\mL}{{\mathcal L}}
\newcommand{\mO}{{\mathcal O}}
\newcommand{\znums}{{\mathbf Z}}
\newcommand{\cnums}{{\mathbf C}}

\newcommand{\sheafhom}{\operatorname{\mathcal{H}\!\mathit{om}}}

\newcommand{\mP}{{\mathbf P}}

\newcommand{\mM}{{\mathcal M}}

\newcommand{\mI}{{\mathcal I}}

\newcommand{\mS}{{\mathcal S}}

\newcommand{\mQ}{{\mathcal Q}}

\newcommand{\Ker}{\operatorname{Ker}}
\newcommand{\mG}{{\mathcal G}}
\newcommand{\bI}{{\mathbf I}}
\newcommand{\mJ}{{\mathcal J}}
\newcommand{\mT}{{\mathcal T}}
\newcommand{\bG}{{\mathbf G}}
\newcommand{\mU}{{\mathcal U}}

\begin{document}

\title[Restricting Bundles to Conics]{Restricting Semistable
Bundles \\ On The Projective Plane To Conics}
\author{Al Vitter}
\address{Department of Mathematics\\
   Tulane University \\
   New Orleans, La. 70118\\
   USA}
\email{vitter@math.tulane.edu}
\thanks{}

\date{Feb. 11, 2004}


\maketitle
\begin{abstract}
We study the restrictions of rank 2 semistable vector bundles $E$ on $\PS{2}$
to conics. A Grauert-M\"ulich type theorem on the generic splitting
is proven. The jumping conics are shown to have the scheme
structure of a hypersurface $J_{2} \subset \PS{5}$ of degree $c_{2}(E)$
when $c_{1}(E)=0$ and of degree $c_{2}(E)-1$ when $c_{1}(E)=-1.$
Some examples of jumping conics and jumping lines are studied in
detail.
\end{abstract}

\section{Introduction} \label{S:intro}

A standard method in the theory of vector bundles on projective
spaces is to restrict a bundle to a line where it splits into
a sum of line bundles and then study how the splitting changes as
the line varies. Grothendieck's result (\cite{splitting},\cite[page 22]{OSS}) gives the
splitting of a rank $r$ bundle $E$ on $\PS{n}$ restricted to a
line $L \subset \PS{n}$ as $E_{L} \cong \bigoplus_{i=1}^{r} \mO_{L}(a_{i})$
for integers $a_{1}\geq a_{2}\dots \geq a_{r}.$ The minimal
(lexicographic order) splitting occurs for a Zariski open subset
(of the Grassmannian $\mG$) of lines in $\PS{n}$ and such lines
are called generic lines for $E.$ The lines on which $E$ has a
larger splitting are called jumping lines and form a proper closed
subscheme $J_{1}$ of $\mG.$ When $E$ is semistable, the
Grauert-Mulich theorem states that the generic splitting satisfies
$a_{j}-a_{j+ 1}\leq 1$ $\forall j$
(\cite[page 192]{OSS} and \cite[in the $r=2$ case]{Barth-sb}).
Further restricting to rank 2
and normalizing $E$ so that $c_{1}(E)=0 \text{ or }-1$, we have that
the generic splitting is $\mO_{L} \oplus \mO_{L}$ when $c_{1}(E)=0$
and $\mO_{L} \oplus \mO_{L}(-1)$ when $c_{1}(E)=-1.$ In the $c_{1}(E)=0$
case, Barth proved \cite{Barth-sb} that the jumping lines form a hypersurface
in $\mG$ of degree $c_{2}(E).$ When $c_{1}(E)=-1$ and $n=2$, Hulek
showed \cite{Hulek} that there are usually only $\binom{c_{2}(E)}{2}$ jumping
lines but that the jumping lines of the second kind form a curve
$\tilde{J_{1}} \subset {\PS{2}}^{*}$ of degree $2(c_{2}(E)-1).$ A
jumping line of the second kind is a line $L$ such that, for $L^{(1)}$
the first order neighborhood of $L$ in $\PS{2}$, $E_{L^{(1)}}$ has
a non-trivial global section.

In this paper we analyze the restriction of rank 2 semistable
bundles $E$ on $\PS{2}$ to conics. One complication arises from
the fact that there are three kinds of conics: smooth conics $C$
isomorphic to $\PS{1}$ via a quadratic parameterization $\PS{1}
\overset{g}{\to}
C$, unions of distinct lines $L_{1}+L_{2}$, and double lines $2L$
(=first order neighborhood $L^{(1)}$ of $L$ in $\PS{2}$). For
smooth conics we prove a Grauert-Mulich type theorem for
semistable $E$ of arbitrary rank (Theorem~\ref{Thm:Grauert-Mulich}). When the rank of $E$
is 2,

$$
c_{1}(g^{*}E)=
\begin{cases}
0 &\text{if }c_{1}(E)=0 \\
-2 &\text{if }c_{1}(E)=-1
\end{cases}
$$

\noindent and so

\begin{equation} \label{E:intro1}
E_{C} \cong g^{*}E \cong
\begin{cases}
\Op{1}(a)\oplus \Op{1}(-a) &\text{for } c_{1}(E)=0 \\
\Op{1}(-1+a)\oplus \Op{1}(-1-a) &\text{for } c_{1}(E)=-1
\end{cases}
\end{equation}

\noindent with $a=0$ corresponding to a generic smooth conic and $a>0$
defining a smooth jumping conic. For the reduced case, jumping
conics are defined as follows. When $c_{1}(E)=0$, $C=L_{1}+L_{2}$ is
a jumping conic means that $L_{1}$ or $L_{2}$ is a jumping line.
When $c_{1}(E)=-1$, $C=L_{1}+L_{2}$ ($L_{1}$ and $L_{2}$ distinct) is
a jumping conic if either $L_{1}$ or $L_{2}$ is jumping or if both
have generic splitting $\Op{1}\oplus \Op{1}(-1)$ and their $\Op{1}$
factors match up at the intersection point $L_{1}\cap L_{2}.$ When $c_{1}(E)=-1$,
and $C=2L$, a jumping conic is a jumping line of the second kind.
We prove that these definitions fit together by showing that the
jumping conics have the subscheme structure of a hypersurface $J_{2}$
in $\PS{5}$, of degree $c_{2}(E)$ if $c_{1}(E)=0$ and of degree
$c_{2}(E)-1$ if $c_{1}(E)=-1$(Theorem~\ref{Thm:jumping conics}). Furthermore, the
singular jumping conics are in the scheme-theoretic closure of the
smooth jumping conics.

The proof of this theorem, easy in the $c_{1}(E)=-1$ case,
requires two substantial ingredients when $c_{1}(E)=0$: 1.the
theory of stable bundles on ruled surfaces and 2.Hurtubise's
analysis of families of rank 2 bundles on $\PS{1}.$ The proof also
gives a new demonstration and explanation of another result of
Barth. Denote the intersection index of a line $\Lambda$ in $\PS{5}$
and $J_{2}$ at $C \in J_{2}$ by $(\Lambda \cdot J_{2})_{C}.$
Define the jump size of $C$ by the integer $a$ in
(\ref{E:intro1}). Then $(\Lambda \cdot J_{2})_{C}\geq a$
(Corollary~\ref{Cor:jump size inequality}).
Actually, Barth's result is the corresponding statement for
jumping lines \cite[Section 6]{Barth-sb}.

Theorem~\ref{Thm:jumping conics} also relates to Hulek's result on the curve $\tilde{J_{1}}$
of jumping lines of the second kind. The double lines comprise the
degree 2 image $\Delta$ of $\mP^{2*}$ in the space of conics $\PS{5}$
via the Veronese map. Since $\tilde{J_{1}}=\Delta \cap J_{2}$ and
since $\text{deg }J_{2}=c_{2}(E)-1$, $\text{deg }\tilde{J_{1}}=2(c_{2}(E)-1).$
Of course this depends on the fact that not all lines are jumping
lines of the second kind and this is the major part of Hulek's
proof.

In Section~\ref{S:(0,2)}, Section~\ref{S:(0,3)}, and Section~\ref{S:(-1,2)},
we describe some examples of jumping conics and jumping lines on rank 2 stable bundles
$E$ on $\PS{2}$ in some detail. When $c_{1}=0$ and $c_{2}=2$, $E$ is
determined by a map into the grassmannian of lines in $\PS{3}$ and
and the jumping lines and conics are described relative to this
map and the Schubert cycles in the grassmannian.
The $c_{1}=0$, $c_{2}=3$ bundles are of two types. The
generic bundle is determined by a map $f \colon \PS{2} \to \PS{2}$
and the jumping lines and conics can be studied in terms of this
map and its ramification divisor. A non-generic bundle is
obtained from $T\PS{2}$ by an elementary modification along a
unique line and its jumping lines and conics can be understood
via this description. Every $c_{1}=-1$, $c_{2}=2$ bundle is an
elementary modification of the trivial bundle along a line. From this we
obtain an explicit equation for its hyperplane of jumping conics.

Ran has obtained \cite{Ran-jumping} a Grauert-Mulich type
theorem and results on jumping curves for rank r
semistable bundles on $\PS{n}$ and their restrictions to rational
curves under certain conditions. His methods are related to
quantum K-theory.

After this paper was written, we learned of the earlier paper of
Manaresi \cite{Manaresi} on the same subject. Our definition of jumping conics is
the same as her definition as is our result on the degree of
the hypersurface of jumping conics. The methods of our paper are
quite different from those of \cite{Manaresi}. In the $c_{1}=0$ case, she constructs a
complete simultaneous deformation of a rank two bundle on two
lines meeting at a point, trivial on one line and non-trivial on
the other, to a trivial bundle on a smooth conic. This is used to
show that the jumping conics form a hypersurface. We do not
understand how her results show that the degree of the
hypersurface is $c_{2}.$

\section{Preliminaries}    \label{S:pre}
By a stable bundle we shall mean Mumford-stable (or $\mu$-
stable), that is

\begin{Def}
Let $X$ be a smooth projective variety of dimension n, $\mO_{X}(1)$
a very ample line bundle on $X$,
and $H$ a corresponding hyperplane section of $X.$ A coherent torsion-free rank r sheaf E
on $X$ is called stable ( resp.\... semistable) if,
for any subsheaf $F\subset E$ of rank
$r_{1}<r$, $c_{1}(F)\cdot H^{n-1}/r_{1}< c_{1}(E)\cdot H^{n-1}/r$
\quad (resp.\... $\leq$ ).
\end{Def}

\noindent When $E$ is a rank 2 bundle on $\PS{2}$ normalized so
that $c_{1}(E)= 0 \text{ or } -1$, stability is equivalent to $h^{0}(\PS{2};E)=0$
when $c_{1}(E)=0$ and to $h^{0}(\PS{2};E(-1))=0$
when $c_{1}(E)=-1.$ For $c_{1}(E)=0$, semistability is equivalent to
$h^{0}(\PS{2};E(-1))=0$ \cite[Ch.2, Sec.1.2]{OSS}. The stable bundles on $\PS{2}$
of fixed Chern classes are
parameterized by a coarse moduli scheme $\mM(c_{1},c_{2})$ which
is a quasi-projective variety \cite{Maruyama2}. By deformation
theory, $\mM(c_{1},c_{2})$ is smooth and $\text{dim }\mM(0,c_{2})=4c_{2}-3$
, $\text{dim }\mM(-1,c_{2})=4c_{2}-4.$

The Riemann-Roch formula for a rank r bundle on $\PS{2}$
\cite[Append.A,Sec.4]{Hartshorne} is

\begin{equation} \label{E:pre1}
\chi(\PS{2};E(k))=
\begin{cases}
\frac{r}{2}(k+2)(k+1)-c_{2}(E) &\text{for } c_{1}(E)=0  \\
\frac{r}{2}(k+2)(k+1)-k-1-c_{2}(E) &\text{for } c_{1}(E)=-1.
\end{cases}
\end{equation}

For a rank r bundle on a smooth genus $g$ curve $Y$, Riemann-Roch
is

\begin{equation} \label{E:pre2}
\chi(Y;E(k))=r(1-g)+rk+c_{1}(E).
\end{equation}

For $E$ a rank 2 bundle on a compact surface $X$, $D \overset{j}{\hookrightarrow}X$
an effective divisor, $\mL$ a line bundle on $D$, and
$\psi \colon E \to E \res{D} \to \mL$ a bundle surjection, define the
rank two bundle $E^{\prime}$ by the elementary modification

\begin{equation} \label{E:pre3}
\begin{CD}
0 @>>> E^{\prime} @>>> E @> \psi >> j_{*}\mL @>>> 0.
\end{CD}
\end{equation}

\noindent Then \cite[page 41]{Friedman}

\begin{align} \label{E:pre4}
c_{1}(E^{\prime}) &=c_{1}(E)-[D]^{*} \quad \in H^{2}(X;\znums)  \\
c_{2}(E^{\prime}) &=c_{2}(E)-[D]^{*}\cdot c_{1}(E)+ c_{1}(\mL)
\quad \in \znums
\end{align}

\noindent where $[D]^{*}$ denotes the cohomology class dual to the
homology class $[D]$ determined by the divisor $D.$

\section{Generic Splitting on Conics}  \label{S:generic}

A conic $C_{\xi}$ in $\PS{2}=\mP(V)$ is defined by the vanishing
of a non-trivial homogeneous polynomial of degree two,
$\xi = \sum_{i,j=0}^{2}\xi_{ij}x_{i}x_{j}=0.$
We let $\xi$ denote this polynomial or the symmetric 3 by 3 matrix
$\xi = (\xi_{ij})$ or the corresponding element of $S^{2}(V^{*})$
as needed. We also use $\xi$ to represent the corresponding
projective class in $\PS{5} \equiv \mP (S^{2}(V^{*})).$ $C_{\xi}$
is a smooth conic if and only if $\det (\xi) \ne 0$ and so the singular
conics form a degree three hypersurface $S \subset \PS{5}.$ For $\xi \in
S$, $C_{\xi}=L+L^{'}$ where $L$ is a line defined by $l(x)=\sum_{i=0}^{2} l_{i}x_{i}=0$
with a similar equation for $L^{\prime}.$ Then ${\PS{2}}^{*} \times {\PS{2}}^{*} \to S$
defined by $(l,l^{\prime})\mapsto l \cdot l^{\prime}$ induces an
isomorphism of the second symmetric power of ${\PS{2}}^{*}$ and $S.$
The diagonal $\Delta$ corresponds to the double lines $C_{\xi}= 2L.$

A conic in $\PS{n}$ for $n\geq 3$ lies in a 2-plane $\PS{2} \subset \PS{n}.$
If $\mG _{2,n}$ is the grassmannian of linear $\PS{2}$'s in
$\PS{n}$ and $\mS$ is the tautological sub-bundle (rank 3) on
$\mG _{2,n}$, the conics are parameterized by $\mP(S^{2}(\mS^{*}))$,
a $\PS{5}$-bundle on $\mG _{2,n}$ of total dimension $3n-1.$

If $C \subset \PS{n}$ is a smooth conic, it has a quadratic
parameterization $g\colon \PS{1} \to C$ unique up to automorphism
of $\PS{1}.$

\begin{Thm} \label{Thm:Grauert-Mulich}
Let E be a rank r semistable bundle on $\PS{n}$ and let
$E_{C}=\bigoplus_{j=1}^{r}\mO_{\PS{1}}(a_{j})$ be the generic
splitting of E on a smooth conic C in $\PS{n}$, with $a_{1}\geq a_{2}\geq \cdots a_{r}.$
Assume that either $n=2$ and r arbitrary or $r=2$ and n
arbitrary. Then $a_{j}-a_{j+1}\leq 1$ for all j.
\end{Thm}

\begin{proof}
The proof of the corresponding (Grauert-Mulich) theorem for
restrictions to lines \cite[Ch.2, Sec.2.1]{OSS} works in our
case with an additional calculation which follows along with a
sketch of the proof. First consider the $n=2$ case. We use the
incidence variety $\bI=\{(x,\xi)\in \PS{2} \times \PS{5} \mid \xi(x)=0\}$
and the projections

$$
\begin{CD}
\bI @>\pi_{1}>> \PS{5} \\
@V\pi_{0}VV \\
\PS{2}
\end{CD}
$$

\noindent
Suppose $a_{i}-a_{i+1}\geq 2$ for some i. We will obtain a
contradiction to E being semistable by constructing a rank i
reflexive
subsheaf $E^{\prime}$ of E such that
$i^{-1}c_{1}(E^{\prime})\cdot \omega_{0} > r^{-1}c_{1}(E)\cdot \omega_{0}.$
This is done by first constructing a rank i reflexive subsheaf ${\hat{E}}^{\prime}$
of $\pi_{0}^{*}E$ on $\bI$ determined as follows.
Set $U \equiv \{\xi \in \PS{5} \mid C_{\xi} \text{ smooth and }
E_{\xi} \text{ splits generically } \}.$ For every
$\xi \in U$, ${\hat{E}}^{\prime}_{\pi_{1}^{-1}(\xi)} \equiv
\bigoplus_{j=1}^{i} \mO(a_{j}).$ The quotient sheaf
$\hat{E}^{\prime \prime} \equiv \pi_{0}^{*} E / {\hat{E}}^{\prime}$ satisfies, for
each $\xi \in U$, $\hat{E}^{\prime \prime}_{{\pi}_{1}^{-1}(\xi)}
\cong
\bigoplus_{j=i+1}^{n} \mO(a_{j}).$ By the descent lemma
\cite[Chap.2, Lemma 2.1.2]{OSS}, there will be a subsheaf $E^{\prime}$
of $E$ on $\PS{2}$ such that $\pi_{0}^{*}E^{\prime} =\hat{E}^{\prime}$ if
  $h^{0}(\bI, \sheafhom(T_{\bI \mid \PS{2}},
\sheafhom(\hat{E}^{\prime},\hat{E}^{\prime \prime})))=0$ where $T_{\bI \mid \PS{2}}$
is the bundle of vectors tangent to the fibers of $\pi_{0}.$ This $E^{\prime}$
clearly yields the desired contradiction. So we must calculate $T_{\bI \mid \PS{2}}.$

The variety $\bI \subset \PS{2} \times \PS{5}$ is defined as the zero set of
the homogeneous polynomial of bidegree $(2,1)$, $f(x,\xi)=x^{t}\xi x$
so $(v,\eta) \in T\PS{2} \oplus T\PS{5}$ is tangent to $\bI$ at $(x,\xi)$
if and only if $v^{t}\xi x + x^{t}\xi v +x^{t} \eta x =0.$ A
tangent to the fiber of $\pi_{0}$ at $(x,\xi)$ has the form $(0,\eta)$
satisfying $x^{t} \eta x =0.$ The Euler sequence

$$
\begin{CD}
0 @>>> \Op{5} @>>\xi> S^{2}V^{*}\otimes \Op{5}(1) @>>> T\PS{5}
@>>> 0
\end{CD}
$$

\noindent shows that $T\PS{5}_{\xi} \cong S^{2}V^{*}/ \cnums \xi$
so $(T_{\bI \mid \PS{2}})_{x,\xi} \cong
\{ \eta \in S^{2}V^{*}/ \cnums \xi \mid x^{t}\eta x=0 \}.$

Set $\Tilde{C}_{\xi} \equiv \pi_{1}^{-1}(\xi)$ for $\xi \in U$,
the copy of the conic $C_{\xi}$ in the incidence variety. To show
$\sheafhom(T_{\bI \mid \PS{2}},\sheafhom(\hat{E}^{\prime},\hat{E}^{\prime \prime}))$
has no global sections it is enough to show it has no sections
along $\Tilde{C}_{\xi}.$ Let $g \colon \PS{1} \to \Tilde{C}_{\xi} \subset \bI$
be a quadratic parameterization. Pulling back $T_{\bI \mid \PS{2}}$
by $g$ we get over $\PS{1}$, with homogeneous coordinates
$u=(u_{0},u_{1})$,

$$
\begin{CD}
0 @>>> g^{*}T_{\bI \mid \PS{2}} @>>>  S^{2}V^{*}/ \cnums \xi
\otimes
\Op{1} @>B>> \Op{1}(4) @>>> 0
\end{CD}
$$

\noindent where $B(u)(\eta) \equiv b(u)\cdot \eta$ for $b(u)$ a
six-dimensional vector function whose components are homogeneous
polynomials of degree 4. Concluding that
$g^{*}T_{\bI \mid \PS{2}} \cong {\Op{1}(-1)} ^{\oplus 4}$ is therefore
equivalent to showing that the map induced by $B$ on sections,
$S^{2}V^{*}/ \cnums \xi \to H^{0}(\PS{1};\Op{1}(4))$,is an
isomorphism. This must be true or else there is a $\xi^{\prime} \in S^{2}V^{*}$
linearly independent of $\xi$ such that $C_{\xi} \subset C_{\xi^{\prime}}$,
which is absurd. Now we have

$$
\sheafhom(T_{\bI \mid \PS{2}},\sheafhom(\hat{E}^{\prime},\hat{E}^{\prime \prime}))
_{\Tilde{C}_{\xi}} \cong \bigoplus_{j\leq i, j^{\prime}>
i} {\Op{1}(-a_{j}+a_{j^{\prime}}+1)}^{\oplus 4}.
$$

\noindent Since $a_{j} \geq a_{j^{\prime}}+2$ for all $j,j^{\prime}$,
there are no sections. This proves our result for bundles over $\PS{2}.$

For $E$ on $\PS{n}$ of rank 2 and $n\geq 3$, Maruyama proved \cite{Maruyama1}
(see also \cite[Sec. 3.3]{Hartshorne-sb} and
\cite[Sec.3.2]{Hartshorne-srs}) that $E_{P}$ is semistable for the
generic plane P in $\PS{n}.$ Now using the incidence variety

$$
\begin{CD}
\bI @>\pi_{1}>> \mP(S^{2}\mS^{*}) \\
@VV \pi_{0} V \\
\PS{n}
\end{CD}
$$

\noindent the same proof works. The same calculation is done on
$\Tilde{C}_{\xi , P} \equiv \pi_{1}^{-1}(\xi,P)$ for $P \equiv \PS{2} \subset \PS{n}$
such that $E_{P}$ is semistable and $\xi$ is a smooth conic in P
on which $E_{P}$ splits generically. A tangent vector to the fiber
of $\pi_{0}$ now has the form $(\eta, t)$ for $\eta \in S^{2}\mS^{*}_{P} / \cnums \xi$
such that $x^{t} \eta x =0$ and $t \in (T\mG_{2,n})_{P}.$

$$
g^{*}T_{\bI \mid \PS{n}} \cong \Op{1}(-1)^{\oplus 4} \oplus
\Op{1}^{\oplus 3(n-2)}.
$$

\end{proof}

\begin{Cor}    \label{Cor:splitting}
Let $E$ be a rank 2 bundle on $\PS{n}$ normalized so that $c_{1} = 0 \text{ or } -1.$
Then $E$ is semistable if and only if for C a generic smooth conic
in $\PS{n}$ and for $\PS{1} \to C$ a quadratic parameterization,

$$
E_{C} \cong
\begin{cases}
\Op{1}\oplus \Op{1} \text{ if } c_{1}=0 \notag  \\
\Op{1}(-1) \oplus \Op{1}(-1) \text{ if } c_{1}=-1.
\end{cases}
$$

\end{Cor}

\begin{proof}
If $E$ is semistable, the splitting follows from
Theorem~\ref{Thm:Grauert-Mulich}. Conversely, if $E$ has the
indicated generic splitting on conics and $\Op{n}(k) \to E$ is a
non-zero bundle map then restricting to a generic conic implies
that $k\leq 0$ if $c_{1}=0$ and $k\leq -1$ if $c_{1} =-1.$
\end{proof}

In connection with Theorem~\ref{Thm:jumping conics}, recall that
indecomposable rank 2 bundles are plentiful on $\PS{2}$ and on $\PS{3}$
\cite{Barth-P2}, \cite{Hartshorne-sb}
but that the only known examples on $\PS{4}$ are variants of the
Horrocks-Mumford bundle \cite{Hulek-hm1}. Furthermore, there are no known examples
on $\PS{n}$ for $n\geq 5$ and the Hartshorne conjecture predicts
that none exist (at least for $n\geq 7$;see \cite{OSS}).

For $r\geq 3, n\geq 3$ there are semistable rank r bundles on $\PS{n}$
whose restrictions to all planes are not semistable (see
\cite{Hartshorne-rt}), e.g. $T\PS{3}.$

\section{Jumping Conics} \label{S:Jumping Conics}
Let E be a rank 2 semistable bundle on $\PS{2}$ normalized so that
$c_{1}(E) =0 \text{ or } -1.$We want to define
a jumping conic $C$ for $E.$ If $C$ is smooth (equivalently
irreducible), this is clear in the light of
Corollary~\ref{Cor:splitting}: If $c_{1}(E)=0$, $E_{C}=\Op{1}(a)\oplus \Op{1}(-a)$
and, if $c_{1}(E)=-1$, $E_{C}=\Op{1}(a-1)\oplus \Op{1}(-a-1).$ In both cases, we define
$C$ to be a jumping conic for $E$ if $a>0.$ We give a separate
definition for singular (equivalently reducible) conics $C=L_{1}+L_{2}$
when $c_{1}(E)=0$:$C$ is jumping if either $L_{1}$
or $L_{2}$ is a jumping line. In the $c_{1}(E)=-1$ case, for $C$ smooth or singular,
define $C$ to be a jumping conic if $h^{0}(C;E_{C})>0.$ Note that
this agrees with our previous definition when $C$ is smooth. When $C$
is singular, $C$ is jumping exactly when either: i)$L_{1}$
or $L_{2}$ is a jumping line, ii)$C=2L$ and $L$ is a jumping line
of the second kind \cite{Hulek}, or iii)$L_{1}$ and $L_{2}$ are generic
so that $E_{L_{j}}=\Op{1}\oplus \Op{1}(-1)$  $j=1,2$ and the $\Op{1}$
summands coincide at $p=L_{1}\cap L_{2}.$

Denote the set of jumping conics by $J_{2}.$
The virtue of these definitions is demonstrated by showing
that $J_{2}$ can be given a scheme structure in
a natural way.

\begin{Thm} \label{Thm:jumping conics}
The set of jumping conics  $J_{2}$ of a rank 2 semistable bundle $E$ on $\PS{2}$
can be given the scheme structure of a hypersurface
in  $\PS{5}$ of degree $c_{2}(E)$ if $c_{1}(E)=0$
and of degree $c_{2}(E)-1$ if $c_{1}(E)=-1.$ Furthermore the
singular jumping conics are in the scheme-theoretic closure of the
smooth jumping conics.
\end{Thm}

\begin{proof}
(The case $c_{1}(E)=-1.$) For $C \subset \PS{2}$
any conic, Riemann-Roch and the sequence

\begin{equation} \label{E:e3.1}
\sexs{E(k-2)}{E(k)}{E_{C}(k)}
\end{equation}

\noindent imply $\chi (C;E_{C}(k))= \chi (\PS{2};E(k))- \chi (\PS{2};E(k-2))=4k.$
This shows that $\pi_{0}^{*}E$ is flat over $\PS{5}$ and that
$h^{0}(C;E_{C})=h^{1}(C;E_{C}).$ Therefore $C$ is a jumping conic if and only if
$h^{1}(C;E_{C})>0.$

The cohomology sequence of (\ref{E:e3.1}) for $k=0$ gives

\begin{equation} \label{E:e3.2}
0 \to H^{0}(C;E_{C}) \to H^{1}(\PS{2};E(-2)) \overset{f_{C}}{\to} H^{1}(\PS{2};E)
\to H^{1}(C;E_{C}) \to 0.
\end{equation}

\noindent From Riemann-Roch and stability,
$h^{1}(\PS{2};E(-2))=h^{1}(\PS{2};E)=c_{2}(E)-1.$
Therefore $C$ is a jumping conic if and only if $\det f_{C}=0.$

Now consider the incidence diagram

$$
\begin{CD}
\bI @>\pi_{1}>> \PS{5} \\
@V\pi_{0}VV \\
\PS{2}
\end{CD}
$$

\noindent which is the restriction to $\bI$ of

$$
\begin{CD}
\PS{2} \times \PS{5} @>p_{1}>> \PS{5} \\
@Vp_{0}VV \\
\PS{2}
\end{CD}
$$

\noindent Note that $E_{0}\equiv \pi_{1*}
\pi_{0}^{*}E =0.$
Set $E_{1}\equiv R^{1}_{\pi_{1}*}\pi_{0}^{*}E.$ Taking the direct image via $p_{1}$
of the exact sequence on $\PS{2} \times \PS{5}$

$$
\begin{CD}
0 @>>> p_{0}^{*}E(-2,-1) @>f>> p_{0}^{*}E @>>> {\pi} _{0}^{*}E
@>>> 0
\end{CD}
$$

\medskip
\noindent gives on $\PS{5}$

$$
\begin{CD}
0 @>>> H^{1}(\PS{2};E(-2))\otimes \Op{5}(-1) @>f_{1}>> H^{1}(\PS{2};E)\otimes \Op{5}
@>>> E_{1} @>>> 0.
\end{CD}
$$

\noindent It follows that
$J_{2}= spt E_{1}=V(\det f_{1})$ as sets and so $\det f_{1}=0$
defines a scheme structure on $J_{2}$, that of a hypersurface in $\PS{5}$
of degree $c_{2}-1.$

Finally, we verify that the smooth jumping conics are dense in $J_{2}.$
The only way this could fail is if the cubic hypersurface of
reducible conics $S \subset \PS{5}$, an irreducible variety, is an irreducible
component of $J_{2}$ (we have seen above that $J_{2}$ has pure
dimension 4, i.e. has no lower dimension components). But this
would contradict Hulek's result \cite{Hulek} that not all lines are jumping
lines of the second kind.
\end{proof}

If $c_{1}(E)=0$, the structure of $J_{2}$ cannot be defined as above
because for $C$ smooth and $E_{C}\cong \Op{1}(a) \oplus \Op{1}(-a)$, $C$ is jumping
iff $a>0$ but $h^{1}(C;E_{C})>0$ iff $a\geq 2.$ Furthermore,
taking direct images on $\PS{5}$ as above,

\begin{multline}
0 \longrightarrow H^{0}(\PS{2};E)\otimes \Op{5} \longrightarrow
E_{0} \longrightarrow H^{1}(\PS{2};E(-2))\otimes \Op{5}(-1) \notag \\
\longrightarrow  H^{1}(\PS{2};E)\otimes \Op{5} \longrightarrow
 E_{1} \longrightarrow  0.
\end{multline}

\noindent where $E_{0}$ is a rank 2 reflexive sheaf,
$h^{1}(\PS{2};E(-2))=c_{2}(E)$, and $h^{1}(\PS{2};E)=c_{2}(E)-2+h^{0}(\PS{2};E).$
Therefore $c(E_{0})= (1-\omega_{1})^{c_{2}(E)} c(E_{1})$ for $\omega_{1}$
the positive generator of $H^{2}(\PS{5};\znums).$ Assume the jump size $a$ is $1$ for all
smooth $\xi \in J_{2}.$ Further assume that the jump size is also $1$
for the reducible jumping conics; we define this to mean that
$h^{0}(L_{1}+L_{2};E_{L_{1}+L_{2}})=2.$ It follows that
$E_{1}=0$ \cite[Ch.3, Cor.12.9]{Hartshorne}, $E_{0}$ is locally free, and
$c(E_{0})=
(1-\omega_{1})^{c_{2}(E)}.$ In particular, $c_{3}(E_{0})=-\binom{c_{2}(E)}{3}.$
This must be zero since $E_{0}$ is locally free of rank 2
and so $c_{2}(E)=0, 1, \text{ or } 2.$
Therefore

\begin{Prop} \label{Prop:jump size two}
Let $E$ be a rank 2 semistable bundle on $\PS{2}$ with $c_{1}(E)=0$
and $c_{2}(E)\geq 3.$ Then E has jumping conics of jump size $\geq 2$
and the support of $R^{1}_{\pi_{1}*}\pi_{0}^{*}E$ has codimension $\leq 3.$
\end{Prop}

To study jumping conics in the $c_{1}(E)=0$ case we need two
tools: 1.Hurtubise's local analysis of families of rank 2 bundles
on $\PS{1}$ \cite{Hurtubise-jumpinglines} and 2.Some results on rank two bundles on
ruled surfaces \cite[Chap.6]{Friedman}.

Let $U$ be an open subset of $\cnums^{n}$ containing the origin
and let $x=(x_{1}, x_{2},\dots x_{n})$ be the coordinates. Let $z$
be the standard affine coordinate on $\cnums \subset \PS{1}.$ For $E$
a rank 2 holomorphic bundle on $\PS{1}\times U$ and for $U$ small
enough, $E$ is trivial on $U_{0}\equiv \{(z,x) \mid z \ne \infty \}$ and
on $U_{1}\equiv \{(z,x) \mid z \ne 0 \}$ and therefore
determined by a 2 by 2 holomorphic matrix transition function $f(z,x)$
on $U_{0}\cap U_{1}.$ Denote by $E_{x}$ the bundle restricted to $\PS{1} \times \{x\}.$

\begin{Prop} \label{Prop:Hurtubise}
(\cite[Prop. 2.1 and 2.4]{Hurtubise-jumpinglines}) Let
$E_{0} \cong \mO(k_{0}) \oplus \mO(-k_{0})$
and let $k\geq k_{0}\geq 0.$ For $U$ chosen small enough, $E$ has a
transition matrix of the form

$$
f(z,x)=\left( \begin{matrix}
z^{k}&p(z,x) \\
0&z^{-k}
\end{matrix} \right)
$$

\noindent for $p(z,x)=\sum_{j=-k+1}^{k-1} a_{j}(x) z^{j}$ [$k=0$ implies $p\equiv 0$]
and $p(z,0)=0$ if $k=k_{0}.$ Define

$$
\Gamma_{p}(x) \equiv \left( \begin{matrix}
a_{0}(x)&a_{1}(x)&\cdots &a_{k-1}(x) \\
a_{-1}(x)&a_{0}(x)&\cdots &a_{k-2}(x) \\
\hdotsfor{4} \\
a_{-k+1}(x)&\cdots &\cdots &a_{0}(x)
\end{matrix} \right)
$$

\noindent Then $E_{x}\cong \mO (a) \oplus \mO (-a)$ iff $\text{rank }\Gamma_{p}(x)=k-a.$
\end{Prop}

Think of $E$ as a family of rank 2 bundles on $\PS{1}$
parameterized by $x \in U.$ If $E$ is trivial on the generic
$\PS{1}$, the jumping lines can then be given a scheme structure
as a hypersurface in $U$, $J=\{x \in U \mid det\Gamma_{p}(x)=0\}.$

Barth \cite[Sect.6]{Barth-sb}gives a different definition of the
jumping line scheme. Applying this to the situation above, take a
resolution of $E(-1)$ on $\PS{1} \times U$ of the form

$$
\sexs{F_{1}}{F_{0}}{E(-1)}
$$

\noindent where $F_{0}=\bigoplus_{i=1}^{r+2}
\pi_{1}^{*}\Op{1}(k_{i})$ for $k_{i}<0$ for each $i$, $\pi_{j}, j=1,2$
the projections from $\PS{1} \times U$, and $F_{1}$
locally free of rank r. Taking direct images on $U$,

$$
\begin{CD}
0 @>>> R^{1}_{\pi_{2}*}F_{1} @>\lambda>> R^{1}_{\pi_{2}*}F_{0}
@>>> R^{1}_{\pi_{2}*}E(-1)  @>>> 0.
\end{CD}
$$

\noindent with the first two terms being bundles of the same rank
and the support of $R^{1}_{\pi_{2*}}E(-1)$ being the jumping locus $J.$ Thus the
scheme structure of $J$ can be defined by $det\lambda =0$ and
Barth shows that this is independent of the resolution.
Furthermore, Hurtubise proves that these two definitions of the
scheme structure agree - this also proves that the Hurtubise
definition is independent of the trivializations used and gives a
second proof that the Barth definition is independent of the
resolution.

\begin{proof}
(The $c_{1}(E)=0$ case.)
We now define the hypersurface structure of $J_{2}$ when $c_{1}(E)=0.$
It is enough to do this locally in $\PS{5}.$ For
$\xi \in \PS{5} \smallsetminus S$, i.e. for $C_{\xi}$ a smooth
conic, use the Hurtubise definition above in a neighborhood of $\xi.$
This defines the jumping locus $J_{2}^{\prime} \subset \PS{5} \smallsetminus S$
as a complex analytic scheme. We will show that $J_{2}^{\prime}$ is
quasi-projective by showing that, as a set, $J_{2}^{\prime}$ is contained
in a closed hypersurface in $\PS{5}.$ Therefore the scheme-theoretic closure of
$J_{2}^{\prime}$ defines a closed sub-scheme of $\PS{5}.$
$J_{2}$ and this closure are shown to be equal as sets, thereby
defining a scheme structure on $J_{2}.$

Consider the second symmetric power of $E$, $S^{2}E$, and note that
$c_{1}(S^{2}E)=0$ and $c_{2}(S^{2}E)=4c_{2}(E).$ For $C$ any
conic,

\begin{multline}  \label{E:e3.3}
0 \longrightarrow H^{0}(\PS{2};S^{2}E) \longrightarrow H^{0}(C;S^{2}E_{C})
\longrightarrow  \\
H^{1}(\PS{2};S^{2}E(-2)) \overset{\alpha_{C}}{\longrightarrow} H^{1}(\PS{2};S^{2}E)
\longrightarrow H^{1}(C;S^{2}E_{C}) \longrightarrow 0.
\end{multline}

\noindent since $H^{0}(\PS{2};S^{2}E(-2))=0$ and
$H^{2}(\PS{2};S^{2}E(-2)) \cong H^{0}(\PS{2};S^{2}E(-1))^{*}=0$ follow
from the semistability of $E.$ From Riemann-Roch,
$h^{1}(\PS{2};S^{2}E(-2))= 4 c_{2}(E)$,
$h^{1}(\PS{2};S^{2}E)= 4c_{2}(E)-3+h^{0}(\PS{2};S^{2}E)$,
and $h^{1}(C;S^{2}E_{C})= h^{0}(C;S^{2}E_{C})-3.$

We show that

\begin{equation} \label{E:e3.4}
C \text{ is a jumping conic iff } h^{1}(C;S^{2}E_{C})>0.
\end{equation}

\medskip

\noindent If $C$ is smooth, $E_{C} \equiv \Op{1}(a)\oplus \Op{1}(-a)$
implies that $h^{0}(C;S^{2}E_{C})$ equals $2a+2$ if $a\geq 1$ and $3$
if $a=0$ and therefore $h^{1}(C;S^{2}E_{C})$ equals $2a-1$ if $a \geq 1$
and $0$ if $a=0.$ For a singular conic $C=L_{1}+L_{2}$, $h^{0}(C;S^{2}E_{C}) \geq 4$
iff $L_{1}$ or $L_{2}$ is a jumping line is also checked by direct
calculation.

Pulling back $S^{2}E$ to $\PS{2}\times \PS{5}$ and to $\bI$ we
have

$$
\sexs{p_{0}^{*}S^{2}E(-2,-1)}{p_{0}^{*}S^{2}E}{\pi
_{0}^{*}S^{2}E}.
$$

\medskip
\noindent Taking direct images by $p_{1}$ gives on $\PS{5}$

\begin{multline} \label{E:e3.5}
0 \longrightarrow H^{0}(\PS{2};S^{2}E)\otimes \Op{5}
\longrightarrow \pi_{1*}\pi_{0}^{*}S^{2}E \longrightarrow \\
H^{1}(\PS{2};S^{2}E(-2))\otimes \Op{5}(-1)
\overset{\alpha}{\longrightarrow} H^{1}(\PS{2};S^{2}E)\otimes \Op{5}
\longrightarrow R^{1}_{\pi_{1}*}\pi_{0}^{*}S^{2}E \longrightarrow
0.
\end{multline}

Let $\mJ$ denote the sheaf $R^{1}_{\pi_{1}*}\pi_{0}^{*}S^{2}E$ and
set $Y= spt \mJ.$
Because $\mJ$ is the cokernel of the bundle homomorphism $\alpha$
and because of (\ref{E:e3.3}) and (\ref{E:e3.4}),
$\xi \in Y$ iff $\alpha_{\xi}$ is not surjective iff
$h^{1}(C_{\xi};S^{2}E_{C_{\xi}})>0$
iff $\xi \in J_{2}.$ Therefore $J_{2}=Y$ as sets and so $J_{2}^{\prime}$
is quasi-projective. We define the scheme structure of $J_{2}$ to
be the closure of $J_{2}^{\prime}$ rather than the scheme
structure of $Y$ for two reasons: 1.The definition of $J_{2}^{\prime}$ via the
method of Hurtubise is more geometrically appealing and 2.The
Hurtubise method allows us to calculate the degree of $J_{2}$ and
to show that the singular jumping conics are in the closure of the
smooth ones.

We will show that the degree of the hypersurface $J_{2} \subset \PS{5}$
is $c_{2}(E)$ by proving that, for the generic line $\Lambda$,
$\Lambda \cdot J_{2}=c_{2}(E).$ Recall that the reducible
(non-smooth) conics form a cubic hypersurface $S\subset \PS{5}$
and that the singular jumping conics form a proper subset of
dimension 3. Therefore the generic line $\Lambda$ misses the
singular jumping conics and so $E$ is trivial on the three
reducible conics of $\Lambda \cong \PS{1}.$ Let $t$ be an affine
parameter on $\Lambda$ and set $X \equiv \pi_{1}^{-1} \Lambda.$
Then $X \to \Lambda$ is a ruled surface with three non-smooth
fibers $C_{t_{i}}$ $i=1,2,3.$ Another way to view $X$ is as
follows. $\Lambda$ is a pencil of conics in $\PS{2}.$ Let $C_{0}$
and $C_{\infty}$ be two smooth conics  with
$C_{0} \cap C_{\infty} = \{p_{0}, p_{1}, p_{2}, p_{3}\}$ the base
locus of $\Lambda.$ $X$ is $\PS{2}$ blown up at these four points and
$C_{t_{i}}$ is the proper transform of $\overline{p_{0}p_{i}} + \overline{p_{j}p_{k}}$
for $(i,j,k)$ a cyclic permutation of $(1,2,3).$ This description
implies that $\pi_{0}^{*}E \res{X}$ has $c_{1}=0$ and $c_{2}=c_{2}(E).$
For simplicity, denote $\pi_{0}^{*}E \res{X}$ by $E.$ Note that the ruled
surface $X$ can be taken to be smooth because a local calculation
shows that smoothness is equivalent to the pencil $\Lambda$ having
four distinct base points.

For all but a finite number of $t \in \Lambda$,
$E_{C_{t}}=\overset{2}{\oplus} \mO_{C_{t}}.$
For $t \in \Lambda \cap J_{2}$, $E_{C_{t}} = \Op{1}(a) \oplus \Op{1}(-a)$
for $a \in \znums^{+}.$ Define the rank 2 bundle $E^{\prime}$ on $X$
as the elementary modification

$$
\begin{CD}
0 @>>> E^{\prime} @>>> E @>\psi>> j_{t*}\Op{1}(-a) @>>> 0
\end{CD}
$$

\noindent for $j_{t}\colon \PS{1} \cong C_{t} \hookrightarrow X$
and $\psi(\sigma) \equiv \Op{1}(-a) \text{-component of } \sigma \res{C_{t}}.$
By (\ref{E:pre4}), $c_{1}(E^{\prime})= -[F]^{*}$
and $c_{2}(E^{\prime})=c_{2}(E)-a$ (where $F$ is the fiber of $X \to \Lambda$).

We use some of the basic theory of rank 2 bundles over ruled
surfaces as developed in \cite[Ch.6]{Friedman}. Some of the
results we quote are proven for geometrically ruled surfaces (all
fibers smooth) but remain true in our case. For $E^{\prime}$ as
above, $E_{C_{t}}^{\prime} = \Op{1}(a^{\prime}) \oplus \Op{1}(-a^{\prime})$
for $0\leq a^{\prime} \leq a$ \cite[page 151]{Friedman}. For any $E$
trivial on the generic fiber of $\pi_{1} \colon X \to \Lambda$, $c_{2}(E) \geq 0$
with equality iff $E$ is trivial on every fiber, i.e. $E=\pi_{1}^{*}W$
for $W$ a rank 2 bundle on $\Lambda$ \cite[Ch.6,Thm.10]{Friedman}.
Therefore after finitely many elementary
modifications of the type described above (more than one
modification may have to be done on the same fiber), we get

$$
\sexs{E^{\prime}}{E}{\mT}
$$

\noindent for $E^{\prime}=\pi_{1}^{*}W$, $W$ a rank 2 bundle on
$\Lambda$, and $\mT$ a torsion sheaf with support on the jumping
fibers. In addition, $c_{2}(E)=\sum_{j} a_{j}$ where the $a_{j} \in {\znums}^{+}$
come from the modifications.

Denote the intersection of $\Lambda$ and $J_{2}$ at $t$ by $(\Lambda \cdot J_{2})_{t}.$
We show that $(\Lambda \cdot J_{2})_{t}=$ the sum of the $a_{j}$'s
coming from the elementary modifications on the fiber $C_{t}.$
This will prove $\Lambda \cdot J_{2}=c_{2}(E).$ For simplicity set
$t=0$ and set $U$ be a small disc about $0$ in $\PS{1}$ with
coordinate $x.$ If $E_{C_{0}}\cong \Op{1}(a) \oplus \Op{1}(-a)$,
apply Proposition~\ref{Prop:Hurtubise} with $k=k_{0}=a$ so that $E$
on $\pi_{1}^{-1}(U)$ is described by the transition function

$$
f(z,x)=\left( \begin{matrix}
z^{a} & p(z,x) \\
0 & z^{-a}
\end{matrix} \right)
$$

\noindent for $p(z,x)=\sum_{j=-a+1}^{a-1} a_{j}(x) z^{j}$ and $a_{j}(0)=0$
for each j. Write $a_{j}(x)=x b_{j}(x)$ and $p(z,x)=x q(z,x).$ If $\mO \oplus \mO$
is the trivialization of $E$ on $U_{0}=\{(z,x) \mid z\ne \infty
\}$, the definition of the elementary modification $E^{\prime}$
implies that the corresponding trivialization of $E^{\prime}$ is
$x\mO \oplus \mO  \subset \mO \oplus \mO.$ The transition function
for $E^{\prime}$ is therefore

$$
\left( \begin{matrix}
x^{-1}&0 \\
0&1
\end{matrix} \right) \cdot \left( \begin{matrix}
z^{a} & p(z,x) \\
0 & z^{-a}
\end{matrix} \right) \cdot \left( \begin{matrix}
x&0  \\
0&1
\end{matrix} \right) = \left( \begin{matrix}
z^{a} & q(z,x) \\
0 & z^{-a}
\end{matrix} \right).
$$

\noindent Note that $\det \Gamma_{p}(x)= x^{a}\det \Gamma_{q}(x).$ If
$E^{\prime}_{C_{0}}$
is trivial, $\det \Gamma_{q}(0)\ne 0$ and $(\Lambda \cdot J_{2})_{t}=a.$
If $E^{\prime}_{C_{0}}\cong \Op{1}(a_{1}) \oplus \Op{1}(-a_{1})$
for $a_{1} \in {\znums}^{+}$, then $a_{1} \leq a$ and we apply
Proposition~\ref{Prop:Hurtubise} again with $k=k_{0}=a_{1}.$ We
get another transition function for $E^{\prime}$ of the form

$$
\left( \begin{matrix}
z^{a_{1}}&r(z,x) \\
0&z^{-a_{1}}
\end{matrix}  \right)
$$

\noindent for $r(z,x)=\sum_{j=-a_{1}+1}^{a_{1}-1} \alpha_{j}(x)z^{j}$
and $\alpha_{j}(0)=0$ for each j. Because Hurtubise's definition
is independent of the trivialization used, $det \Gamma_{q}(x)$ and $det \Gamma_{r}(x)$
have the same order zero at $x=0.$ $det \Gamma_{r}(x)=x^{a_{1}} det \Gamma_{s}(x)$
for $r(z,x)=x s(z,x).$ Therefore the order of the zero of $det \Gamma_{p}(x)$
at $x=0$ is $a+a_{1}+ \text{ order of }det \Gamma_{s}(x)$ at $x=0.$
Our result follows by induction.

Finally, we prove that $\overline{J_{2}^{\prime}}=J_{2}.$ For $\psi
\colon J_{1}\times {\PS{2}}^{*} \to S \subset \PS{5}$ defined by
$\psi(l,l^{\prime})=l\cdot l^{\prime}$, the image of $\psi$ is $J_{2}\cap S.$
Let $J_{1}=\sum_{\mu}m_{\mu}J_{1\mu}$ as a divisor; $J_{1\mu}$ is
an irreducible component of $J_{1}$ for each $\mu.$ Then the
irreducible components of $S\cap J_{2}$ are the 3-dimensional
varieties $Y_{\mu}\equiv \psi(J_{1\mu}\times {\PS{2}}^{*}).$ Suppose
$\eta \in Y_{\mu}$ is not in $\overline{J_{2}^{\prime}}.$ Since
$\psi(J_{1}\times J_{1})$ is 2-dimensional, we can assume
$\eta \notin \psi(J_{1}\times J_{1})$, that is, $C_{\eta}=L_{1}+L_{2}$
for $L_{1}$ a generic line and $L_{2}$ a jumping line. Since the
set of conics tangent to $C_{\eta}$ at some point is
4-dimensional, the generic choice of $\xi \in \PS{5}$ corresponds
to a smooth $E$-generic conic not tangent to $C_{\eta}.$ Letting $\Lambda$
be the line in $\PS{5}$ through $\eta$ and $\xi$, we can also
assume $\Lambda \cap (\overline{J_{2}^{\prime}}\cap S)=\emptyset.$
Let $X=\pi_{1}^{-1}\Lambda$ be the corresponding ruled surface and
denote $\pi^{*}_{0}E\res{X}$ by $E$ as above. The jumping fibers of $X \to \Lambda$
include $C_{\eta}$ and the fibers above the points of
$\overline{J_{2}^{\prime}}\cap \Lambda$, the latter being smooth.
From $C_{\eta}^{2}=0$ and $C_{\eta}\cdot K_{X}=-2$ it follows that
$L_{i}^{2}=-1$ and $L_{i}\cdot K_{X}=-1$ for $i=1,2.$ Blow down
the generic line $L_{1}$ to produce a ruled surface $\Check{X}.$ $L_{2}$
blows down to a smooth fiber $\Check{L}_{2}$ of $\Check{X}.$ By a theorem of
Schwarzenberger \cite [Thm.5]{Schwarzenberger1}, the bundle $E$
descends to a bundle $\Check{E}$ on $\Check{X}$ ($E=\pi^{*}\Check{E}$
for $\pi \colon X \to \Check{X}$ the blow-down map) because $E$ is
trivial $L_{1}.$ Note that
$c_{2}(\Check{E})=c_{2}(E)$, $c_{1}(\Check{E})=c_{1}(E)$, and
$\Check{E}_{\Check{L}_{2}} \cong \mO_{\Check{L}_{2}}(b) \oplus \mO_{\Check{L}_{2}}(-b)$
for $b\in {\znums}^{+}$ (if $b=0, E_{L_{2}}=\pi^{*}\Check{E}_{\Check{L}_{2}}$
would be trivial, a contradiction). Also note that the splitting
of $E$ on $C_{t}$ is the same as that of $\Check{E}$ on $\Check{C}_{t}=\pi(C_{t})$
for all $t\ne \eta.$

Our previous calculation shows that, for each $t \in \Lambda$ such
that $C_{t}$ is smooth, $(\Lambda \cdot \overline{J_{2}^{\prime}})_{t} = \sum a_{j}$
where the sum is over all $a_{j}\in {\znums}^{+}$ coming from the
elementary modifications involving the fiber $C_{t}$ needed to
change $E$ so that it becomes trivial on $C_{t}.$ We used this to
prove that $\Lambda \cdot \overline{J_{2}^{\prime}}= c_{2}(E)$ for
the generic $\Lambda.$ Since the $\Lambda$ chosen above is not
contained in $\overline{J_{2}^{\prime}}$,
$\Lambda \cdot \overline{J_{2}^{\prime}}= c_{2}(E).$ Doing the
same calculation for $\Check{E}$ gives

\begin{align}
c_{2}(E)&=c_{2}(\Check{E})  \notag \\
\quad &\geq b+ \sum_{t\in \Lambda \cap \overline{J_{2}^{\prime}}}
(\Lambda \cdot \overline{J_{2}^{\prime}})_{t}  \notag \\
\quad &=b+\Lambda \cdot \overline{J_{2}^{\prime}} \notag \\
\quad &=b+c_{2}(E)     \notag
\end{align}

\noindent which is a contradiction. Therefore $\overline{J_{2}^{\prime}}=J_{2}.$
\end{proof}

\begin{Cor} \label{Cor:jump size inequality}
Let $E$ be a stable rank 2 bundle on $\PS{2}$ with $c_{1}(E)=0.$
Let $C$ be a smooth jumping conic of $E$ of jump size $a$ corresponding to
$\xi \in J_{2} \subset \PS{5}.$
Let $\Lambda$ be a line in $\PS{5}$ through $\xi$ not contained in
$J_{2}.$ Then the intersection multiplicity of $\Lambda$ and $J_{2}$
at $\xi$ satisfies $(\Lambda \cdot J_{2})_{\xi} \geq a.$
\end{Cor}

\section{Rank 2 Stable Bundles on $\PS{2}$ with $c_{1}=0$,
$c_{2}=2$} \label{S:(0,2)}

Let $E$ be a rank 2 stable bundle on $\PS{2}$ with $c_{1}=0 \text{
and } c_{2}=2.$ Set $W\equiv H^{0}(\PS{2};E(1))$ and note that $\sigma \in W$
has 3 zeros. The sequence

$$
\begin{CD}
0 @>>> \Op{2}(-1) @> \sigma >> E @> \sigma \wedge >>
\mI_{Z_{\sigma}}(1) @>>> 0
\end{CD}
$$

\noindent and the stability of $E$ imply that the zeros of $\sigma$
are not colinear.
Friedman \cite{Friedman} shows that $\text{dim } W=4$ and that $E$ can be
described by a sequence

\begin{equation} \label{E:e(0,2)1}
\begin{CD}
0 @>>> W^{\prime} \otimes \Op{2}(-1) @> f >> W \otimes \Op{2} @> e >> E(1)
@>>>0
\end{CD}
\end{equation}

\noindent where $e$ is evaluation of sections and $f$ is a 4-by-2
matrix of homogeneous degree 1 polynomials.
Denote the two columns
of $f(x)$ by $f_{j}(x)$, $j=1,2.$ We will also use $f$ to denote
the induced map $f \colon \PS{2} \to \bG$
($\bG$ is the grassmannian of 2-planes through $0 \in W$ or the lines in
$\mP W \cong \PS{3}$)
and write $f(x)=f_{1}(x)\wedge f_{2}(x).$ $f(x)=\{ \sigma \in W \mid \sigma (x)=0 \}.$
Note that $E(1) \cong f^{*}\mQ$ for $\mQ$ the standard rank-2 quotient
bundle on $\bG$ so that $f$ determines $E.$ Also notice that $f^{*}\mQ \cong g^{*}\mQ$
iff there is an $A \in GL_{4}$ and a $B \in SL_{2}$ such that $g=AfB.$
The parameter count $dim\{ f \} -dim\{ A \} -dim\{ B \} =8\cdot 3-16-3=5$
agrees with $\text{dim }\mM(0,2).$
We will describe the jumping lines and conics of $E$ in terms of
the map $f.$

First recall the Schubert cycles whose classes generate the
cohomology of $\bG$ \cite[Ch.1,Sec.5]{GH}. Let $\Lambda \subset \mP W$ be
the line corresponding to $l \in \bG$. A line $\Lambda_{0} \subset \mP W$
defines the codimension one cycle
$Z_{1}(\Lambda_{0}) \equiv \{l\in \bG \mid \Lambda_{0} \cap \Lambda \ne \emptyset \}$;
its dual cohomology class is given by the Fubini-Study form $\omega_{1}.$
A point $w \in \mP W$ defines the codimension two cycle
$Z_{2}(w)\equiv \{l\in \bG \mid w \in \Lambda \}$; its dual
class is denoted $\omega_{2}.$ A 2-plane $P \subset \mP W$ defines
the codimension two cycle $Z_{1,1}(P)\equiv \{l\in \bG \mid \Lambda \subset P \}$;
its dual class is denoted $\omega_{1,1}.$ We have

\begin{equation} \label{E:e(0,2)2}
\omega_{1}^{2}=\omega_{2}+\omega_{1,1}
\end{equation}

\noindent and

\begin{equation} \label{E:e(0,2)3}
c(\mQ)=1+\omega_{1}+\omega_{2}.
\end{equation}

It is sometimes useful to regard $W$ as an abstract 4-dimensional vector space
and $E$ as being the quotient bundle defined by (\ref{E:e(0,2)1}).
The cohomology sequence of (\ref{E:e(0,2)1}) then implies that global sections
$\sigma$ have the form $\sigma = \sigma_{w}$ for $w\in W$ and
$\sigma_{w}(x)=w \text{ mod }f(x).$
Therefore $Z_{\sigma_{w}}=f^{-1}(Z_{2}(w))$ and so

\begin{equation} \label{E:e(0,2)4}
\int_{\PS{2}} f^{*}\omega_{2} =3.
\end{equation}

\noindent Since $f$ is quadratic,

\begin{equation} \label{E:e(0,2)5}
\int_{\PS{2}} f^{*}\omega_{1}^{2} =\int_{\PS{2}}
(2\omega_{0})^{2}=4
\end{equation}

\noindent and therefore (\ref{E:e(0,2)2})  and
(\ref{E:e(0,2)4})imply

\begin{equation} \label{E:e(0,2)6}
\int_{\PS{2}} f^{*}\omega_{1,1} =1
\end{equation}

\noindent and $f(\PS{2})\cdot Z_{1,1}(P)=1$ for the generic
2-plane $P.$

Let $L$ be a line in $\PS{2}$ and let $p$ and $q$ be distinct points on $L.$
$E(1)_{L}\cong \mO_{L}(1)\oplus \mO_{L}(1)$ for $L$ generic and
$E(1)_{L}\cong \mO_{L}(2)\oplus \mO_{L}$ for $L$ a jumping line
(any other splitting contradicts $E(1)$ being globally generated).

\begin{Prop}
Let $E$ be a rank 2 stable bundle on $\PS{2}$ with $c_{1}=0$, $c_{2}=2$
and let $f \colon \PS{2} \to \bG$ be the corresponding map. The line $L=\overline{pq}$
is a jumping line for $E$ iff $f(p)\wedge f(q)=0$ iff there is a
Schubert cycle $Z_{1,1}(P)$ in $\bG$ such that $L=f^{-1}(Z_{1,1}(P)).$
\end{Prop}

\begin{proof}
Parameterize $L$ by $t_{0}p+t_{1}q$; then $f$ restricted to $L$ is

$$
f_{L}(t)=(t_{0}f_{1}(p)+t_{1}f_{1}(q))\wedge
(t_{0}f_{2}(p)+t_{1}f_{2}(q)).
$$

We collect some facts about $f.$ $f$ is 1-to-1 because if
$f(p)=f(q)=u \wedge v \in \bG$, $f$ is constant on $\overline{pq}$
and so $\sigma_{u}$ vanishes on $\overline{pq}$ contradicting the
stability of $E.$ The image of $f$ does not lie in
any cycle $Z_{1}(\Lambda_{0})$ nor therefore in any cycle $Z_{1,1}(P).$
This holds because if $f(\PS{2}) \subset Z_{1}(\Lambda_{0})$,
$\PS{2} \subset \bigcup_{w \in \Lambda_{0}}Z_{\sigma_{w}}$ which
is impossible by a dimension count.

\noindent If $f(p)\wedge f(q) \ne 0$, $f_{1}(p)$, $f_{1}(q)$,
$f_{2}(p)$, $f_{2}(q)$ form a basis for $W.$
It follows that $\sigma \in W$ can have at
most one zero on $L$ which implies that $L$ is generic.

If $f(p)\wedge f(q)=0$, the span of $f_{1}(p)$, $f_{1}(q)$,
$f_{2}(p)$, $f_{2}(q)$ has dimension 3 (because $f$ is 1-to-1).
Choose a basis $w_{i}$, $i=0 \text{ to } 3$ for W so that

$$
f_{L}(t)= (t_{0} w_{0}+t_{1} w_{1})\wedge (t_{0} w_{2}+t_{1}[a
w_{0}+b w_{1}+c w_{2}])
$$

\noindent It follows that $\sigma_{w_{2}}$ has two zeros on $L$ if
$b \ne 0$ and $\sigma_{w_{0}}$ has two zeros on $L$ if $b=0$
and so $L$ is a jumping line. Furthermore $L=f^{-1}(Z_{1,1}(P))$
where $P$ is the $\PS{2}$ spanned by $w_{i}$, $i=0,1,2.$
Containment is clear. If $r \notin L$ satisfies $f(r) \in
Z_{1,1}(P)$, then $f(\PS{2}) \subset Z_{1,1}(P)$ which we know is
impossible.
\end{proof}

\begin{Prop}
Let $E$ be a rank 2 stable bundle on $\PS{2}$ with $c_{1}=0$, $c_{2}=2$
and let $f \colon \PS{2} \to \bG$ be the corresponding map. For a
conic $C \subset \PS{2}$, the following are equivalent: \\
i) C is a jumping conic for $E$. \\
ii) There are global sections $\sigma$ and $s$ of $E(1)$ such that
$C=Z_{\sigma \wedge s}.$ \\
iii) There is a Schubert cycle $Z_{1}(\Lambda_{0}) \subset \bG$
such that $C=f^{-1}(Z_{1}(\Lambda_{0})).$
\end{Prop}

\begin{proof}
Let $C$ be a smooth conic and give it a quadratic parameterization
$\PS{1} \to C.$ $E(1)_{C}\cong \Op{1}(2+a)\oplus \Op{1}(2-a)$ and
$C$ is jumping iff $a\geq 1.$ Set
$W_{1}\equiv \{ \sigma \in W \mid \Op{1}(2-a)-\text{component of }\sigma =0 \}.$
If $C$ is a jumping conic, $\text{dim }W_{1}\geq 2.$ For
$\sigma$, $s\in W_{1}$ independent, $C=Z_{\sigma \wedge s}.$
Conversely, if $C=Z_{\sigma \wedge s}$, $\sigma$ and $s$ each have
all three of their zeros on $C$ so $a=1$ and $C$ is a jumping
conic. Recall that $\sigma$ has the form $\sigma_{u}$ for $u \in W$
and $\sigma_{u}(x)=u \text{ mod }f(x)$ and, similarly, $s=s_{v}.$
Thus the homogeneous equation of degree 2 for a jumping conic $C$
can be written $u \wedge v \wedge f(x)=0.$ It follows that
$C=f^{-1}(Z_{1}(u \wedge v)).$

We have shown that the image of the map from $\bG$ to $J_{2}$ given by
$u \wedge v \to Z_{\sigma_{u} \wedge s_{v}}$ contains all the
smooth jumping conics and therefore the map is surjective.
\end{proof}

\section{Rank 2 Stable Bundles on $\PS{2}$ with $c_{1}=0$,
$c_{2}=3$} \label{S:(0,3)}

Let $E$ be a rank $2$ stable bundle on $\PS{2}$ with $c_{1}=0 \text{ and } c_{2}=3.$
From Riemann-Roch, $h^{0}(\PS{2};E(1))\geq 3$ and, for a non-zero
section $\sigma \in H^{0}(\PS{2};E(1))$, the zero set $Z_{\sigma}$
is a 0-dimensional subscheme of length $4$ and

$$
\begin{CD}
0 @>>> \Op{2} @> \sigma >> E(1) @> \sigma \wedge >>
\mI_{Z_{\sigma}}(2) @>>> 0.
\end{CD}
$$

\noindent Now assume that no three points of $Z_{\sigma}$ are
collinear. In this case, we say that $E$ is a bundle of general type
and prove later that all sections of $E$ have zero sets with no
three points collinear.
Friedman \cite[page 94]{Friedman}
points out that $h^{0}(\PS{2};E(1))=3$ and that $E(1)$ is globally
generated, that is, for $W \equiv H^{0}(\PS{2};E(1))$,

\begin{equation} \label{E:e(0,3)1}
\begin{CD}
0 @>>> \Op{2}(-2) @> f >> W \otimes \Op{2} @> e >> E(1) @>>> 0.
\end{CD}
\end{equation}

\noindent where $e$ is evaluation of sections and $f$ is defined
by a 3-vector of degree 2 homogeneous polynomials. We also use $f$ to denote the
induced regular map

\begin{equation} \label{E:e(0,3)2}
f\colon \PS{2} \to \mP W \cong \PS{2},
\end{equation}

\noindent a 4-to-1
branched cover. Note that $f(x)$ equals the point in $\mP W$
corresponding to the 1-dimensional space of sections that vanish
at $x.$ [If $\sigma \in W \setminus 0$ satisfies $\sigma(x) =0$ we
will also write $f(x)=\sigma.$] For $\mQ$ the canonical quotient
bundle on $\PS{2}$, $f^{*}\mQ =E(1)$ so that $f$ determines $E$
and $f^{*}\mQ = g^{*}\mQ $ iff there is an $A \in GL_{3}$ such
that $g=Af.$ Counting parameters, $dim\{E\}= dim\{f\} - dim\{A\} = 3\binom{4}{2} -9 =9$
which matches $\text{dim } \mM (0,3)$ \cite{Barth-P2}. By the Hurwitz formula, the
ramification divisor $R$ of $f$ is a cubic curve.

We will describe the jumping lines and conics of $E$ in terms of the
map $f$ and its ramification divisor $R.$ Define the 3 by 3
matrix of degree 1 homogeneous polynomials in $x_{0}, x_{1},
x_{2}$, $F\equiv ({f_{i}}_{x_{j}})$ and set
$D\equiv det F = f_{x_{0}}\wedge f_{x_{1}} \wedge f_{x_{2}}.$ $R$
is defined by $D(x)=0.$ On an affine piece of $\PS{2}$, for
example on $U_{0}\equiv \{x \in \PS{2} \mid x_{0} \neq 0 \}$, we
use Euler's identity to write $D$ in inhomogeneous form
$D=3 f \wedge f_{x_{1}} \wedge f_{x_{2}}$ where $x_{1}, x_{2}$
are the affine coordinates obtained by setting $x_{0}=1.$

For $p \in \PS{2}$ and $f(p)=\sigma \in \mP W$, make a change in
the homogeneous coordinates in both the domain and range to obtain
$p=e_{0}$ and $e_{0}, e_{1}, e_{2}$ a basis for $W$ with $e_{0}=\sigma.$
The global sections $e_{1}$ and $e_{2}$ form a local frame for $E(1)$
and so $\sigma = \sigma_{1} e_{1}+\sigma_{2} e_{2}$
near $p=e_{0}.$ From

$$
f_{0}e_{0}+f_{1}e_{1}+f_{2}e_{2}=0
$$

\noindent it follows that

$$
\sigma_{i}(x) =-f_{i}(x)/f_{0}(x)
$$

\noindent for $i=1,2$ near $p=e_{0}$ and therefore $df_{p}=-d\sigma_{p}.$
Defining the index of the
zero of $\sigma$ at $p$ to be $i_{p}(\sigma)\equiv dim
\mO_{\PS{2},p}/(\sigma_{1},\sigma_{2})_{p}$, we see that $p \in R$
iff $p$ is a zero of $\sigma =f(p)$ and $i_{p}(\sigma)\geq 2.$

Now suppose $p \in R$ so that $\text{rank }F_{p} \leq 2.$ If $\text{rank }F_{p}=1$,
$R$ is clearly singular at $p$ and $dD_{p}=0$, $df_{p}=0$, and $d\sigma_{p}=0.$
If $\text{rank }F_{p}=2$, further changes in the homogeneous
coordinates allows us to assume $f_{x_{1}}(e_{0})=e_{1}$ and $f_{x_{2}}(e_{0})=0.$
Therefore

\begin{align}  \label{e(0,3)7}
df_{e_{0}}&= \left(\begin{matrix}
1 0 \\
0 0
\end{matrix} \right)   \notag  \\
d\sigma_{e_{0}}&=\left(\begin{matrix}
1 0 \\
0 0
\end{matrix} \right).
\end{align}

\noindent with respect to the obvious bases and

\begin{equation} \label{E:e(0,3)8}
dD_{e_{0}}=3{f_{2}}_{x_{2}x_{1}}(e_{0}) dx_{1}+3{f_{2}}_{x_{2}x_{2}}(e_{0})
dx_{2}.
\end{equation}

\noindent Thus $R$ is singular at $e_{0}$ iff ${f_{2}}_{x_{2}x_{1}}(e_{0})=0$
and ${f_{2}}_{x_{2}x_{2}}(e_{0})=0$. If $R$ is smooth at $e_{0}$,
$TR_{e_{0}}=\Ker df_{e_{0}}=\cnums e_{2}$ iff ${f_{2}}_{x_{2}x_{2}}(e_{0})=0$.

\begin{Prop}
The ramification curve $R$ is either a smooth cubic, the union
of a line and a smooth conic, or the union of three distinct
lines. The last case occurs exactly when there is a point $p \in R$
satisfying $\text{rank }F_{p}=1.$
\end{Prop}

\begin{proof}
Let $p$ be a singular point of $R.$ If $\text{rank }F_{p}=1$, then
choosing coordinates correctly, we get $p=e_{0}$,
$f(e_{0})=e_{0}$, and $f_{x_{1}}(e_{0})=0=f_{x_{2}}(e_{0}).$ Write
$f_{l}=\sum_{0 \leq i \leq j \leq 2} a_{lij}x_{i}x_{j}.$ Imposing the above
conditions we obtain

$$
f(x)=\left( \begin{matrix}
x_{0}^{2}+ a_{001}x_{0}x_{1}+a_{002}x_{0}x_{2}
+\sum_{1 \leq i \leq j \leq 2} a_{0ij}x_{i}x_{j} \\
\sum_{1 \leq i \leq j \leq 2} a_{1ij}x_{i}x_{j} \\
\sum_{1 \leq i \leq j \leq 2} a_{2ij}x_{i}x_{j}
\end{matrix} \right)
$$

\noindent My making an additional linear change in $x_{0}, x_{1}, x_{2}$ we
can put $f$ in one of the following two forms:

$$
f(x)=\left( \begin{matrix}
x_{0}^{2} + \sum_{1 \leq i \leq j \leq 2} a_{0ij}x_{i}x_{j} \\
x_{1}x_{2} \\
(x_{2}-\alpha x_{1})(x_{2}-\beta x_{1})
\end{matrix} \right)
$$

\noindent for $\alpha \neq 0$ and $\beta \neq 0$ (this follows from the fact that $f$
has finite fibers) or

$$
f(x)=\left( \begin{matrix}
x_{0}^{2} + \sum_{1 \leq i \leq j \leq 2} a_{0ij}x_{i}x_{j} \\
x_{1}^{2} \\
x_{2}^{2 }
\end{matrix} \right).
$$

\noindent In the first case,

$$
D=4x_{0}(x_{2}^{2}-\alpha \beta x_{1}^{2})
$$

\noindent and we see that $R$ is the union of three distinct
lines. The second case gives the same conclusion.

If $p$ is a singular point of $R$ and $\text{rank }F_{p}=2$,
choose coordinates as above to obtain

$$
f(x)=\left( \begin{matrix}
x_{0}^{2}+\sum_{1 \leq i \leq j \leq 2} a_{0ij}x_{i}x_{j} \\
x_{0}x_{1}+\sum_{1 \leq i \leq j \leq 2} a_{1ij}x_{i}x_{j} \\
a_{211}x_{1}^{2}
\end{matrix} \right)
$$

\noindent The fact that $a_{212}=0$ and $a_{222}=0$ was shown above.
Note that $a_{211}\ne 0$ and $a_{122}\ne 0$ follow from the fact that $f$
has finite fibers. Direct calculation gives

$$
D=2a_{211}x_{1}[2x_{0}(2a_{122}x_{2} +
a_{112}x_{1})-x_{1}(2a_{022}x_{2} + a_{012}x_{1})]
$$

\noindent and therefore $R=L+C$ for $L$ the line defined by $x_{1}=0$
and $C$ a conic. The conic is degenerate iff
$a_{122}(a_{122}a_{012}-a_{112}a_{022})=0$, i.e. $a_{122}a_{012}-a_{112}a_{022}=0.$
In this case, make the linear change of homogeneous coordinates in the range given by
the matrix

$$
\left( \begin{matrix}
1 & 0 & -a_{011}/a_{211} \\
0 & 1 & -a_{111}/a_{211} \\
0 & 0 & 1
\end{matrix} \right)
$$

\noindent and then another given by

$$
\left( \begin{matrix}
1 & -\lambda & 0 \\
0 & 1 & 0 \\
0 & 0 & 1
\end{matrix} \right)
$$

\noindent for $\lambda$ defined by $\lambda (a_{122}, a_{112})=(a_{022}, a_{012}).$
This yields

$$
f=\left( \begin{matrix}
x_{0}^{2}-\lambda x_{0}x_{1} \\
x_{0}x_{1}+a_{112}x_{1}x_{2}+a_{122}x_{2}^{2} \\
a_{211}x_{1}^{2}
\end{matrix} \right)
$$

\noindent Calculating $F$ and $D$ shows that $R$ is the union of three
distinct lines and that $\text{rank }F_{e_{2}}=1.$
\end{proof}

For any line $L \subset \PS{2}$, $E(1)_{L} \cong \mO_{L}(1+a) \oplus \mO_{L}(1-a)$
with $0 \leq a.$ For $L$ generic, $a=0.$ Since $E(1)$ is globally
generated, $a\leq 1$ so the jumping lines correspond to $a=1.$
Note that this means that for every global section $\sigma$ of
$E(1)$, no 3 of the 4 zeroes of $\sigma$ are collinear.
Let $f_{L}$ be the restriction $f$
to $L$ and denote the image by $\hat{L}.$

If $L$ is generic, $E(1)_{L} \cong \mO_{L}(1) \oplus \mO_{L}(1)$
implies that $f_{L}$ is 1-to-1. Since $f$ is quadratic and $L$ is irreducible,
$\hat{L}$ is a smooth conic.

For $L$ jumping, $E(1)_{L} \cong \mO_{L}(2) \oplus \mO_{L}$
implies that $f_{L}$ is 2-to-1 and that $\hat{L}$ is a line. By
the Hurwitz theorem, $f_{L}$ has 2 ramification points, 2 of the 3
points of $L\cdot R.$ We now work in the opposite direction by
starting with a point of $R$ and producing a jumping line.
For $p \in R$ with $\text{rank } {d\sigma}_{p}=1$,
let $v \in T\PS{2}_{p}$ span $\Ker df_{p}=\Ker {d\sigma}_{p}.$ For $L$ the
line through p with direction $v$ and $\sigma =f(p)$, $\sigma_{L}$
clearly has a double zero at p and so $E(1)_{L} \cong \mO_{L}(2) \oplus \mO_{L}.$
$L$ meets $R$ at 2 more points. One of them, denoted by q, must be distinct from $p$
and satisfy $TL_{q}=\Ker{d\tau}_{q}$
for $\tau =f(q).$ Thus p and q are the ramification points of the
2-to-1 branched cover $f_{L} \colon L \rightarrow \hat{L}.$

\begin{Prop}
Let $E$ be a rank 2 stable bundle on $\PS{2}$ with $c_{1}=0$, $c_{2}=3$
of general type. Let $f \colon \PS{2} \to \PS{2}$ be the
corresponding map with ramification curve $R$. If $R$ is smooth,
there is a regular 2-to-1 map $g \colon R \to J_{1}$ onto the cubic curve
of jumping lines of $E.$
\end{Prop}

\begin{proof}
The definition of the jumping line $g(p)$ corresponding to $p \in R$ was given above.
Because $R$ is smooth,  $\text{rank }df_{p}=1.$ As in
the previous proof, we can change homogeneous coordinates in the
domain and range so that $p=e_{0}$, $f(e_{0})=e_{0}$,
$f_{x_{1}}(e_{0})=e_{1}$, and $f_{x_{2}}(e_{0})=0.$ In local
coordinates $x_{1}, x_{2}$ near $p$,

$$
f(x)=\left( \begin{matrix}
1 \\
f_{1}(1,x_{1},x_{2})/f_{0}(1,x_{1},x_{2}) \\
f_{2}(1,x_{1},x_{2})/f_{0}(1,x_{1},x_{2})
\end{matrix} \right)
$$

\noindent For $x \in R$ near $p$, the kernel of $df_{x}$ is
spanned by

$$
v(x)=\left( \begin{matrix}
0 \\
-(f_{1}/f_{0})_{x_{2}} \\
(f_{1}/f_{0})_{x_{1}}
\end{matrix} \right)
$$

\noindent For $x \in R$ near $p$, $g(x)=x \wedge v(x)$ and so $g$
is a regular mapping.
\end{proof}

Let $C \subset \PS{2}$ be a smooth conic. $E(1)_{C} \cong \Op{1}(2+a) \oplus \Op{1}(2-a)$
for $a \geq 0.$ Because $E$ is globally generated, $a \leq 2.$ The
smooth jumping conics have $a=1 \text{ or } 2.$ If $a=2$ and $\xi \in J_{2}$
is the point corresponding to $C$, $J_{2}$ is singular at $\xi$ by
Corollary~\ref{Cor:jump size inequality}.

Denote the restriction of $f$ to $C$ by $f_{C}$ and its image by $\hat{C}.$

First consider a smooth generic conic: $E(1)_{C} \cong \Op{1}(2) \oplus
\Op{1}(2).$ For $\omega_{0}$ the Fubini-Study form on $\PS{2}$,

$$
\int_{C}f^{*}_{C} \omega_{0} = \int_{C}2\omega_{0} = 4
$$

\noindent so either i) $deg \hat{C} =1$ and $f_{C}$ is
generically 4-to-1, ii) $deg \hat{C} =2$ and $f_{C}$ is generically
2-to-1, or iii) $deg \hat{C} =4$ and $f_{C}$ is generically
1-to-1. Case i) can not occur because a section over $C$ can have
at most 2 zeroes. Case ii) can not occur or else for $(\sigma_{1},\sigma_{2})$
defined as the image of $\sigma$ by

$$
0 \to H^{0}(\PS{2};E(1)) \overset{rest.}{\to} H^{0}(C;E(1)_{C})
\cong \overset{2}{\oplus} H^{0}(\PS{1};\Op{1}(2))
$$

\noindent $\sigma_{1}=c \sigma_{2}$ for $c \in \cnums^{*}$ (c
could depend on $\sigma$). But then, for any $x \in C$,
$dim\{ \sigma \mid \sigma (x)=0 \}=2$, which contradicts $E(1)$
being globally generated. Thus $f_{C} \colon C \to \hat{C}$ is
generically 1-to-1 and $deg \hat{C} =4.$ It follows that the
arithmetic genus of $\hat{C}$ is 3, that $\hat{C}$ is singular,
that $f_{C}$ is a normalization of $\hat{C}$, and that, for

$$
\begin{CD}
0 @>>> \mO_{\hat{C}} @> f^{*}_{C} >> f_{C*}\mO_{C} @>>>
f_{C*}\mO_{C}/\mO_{\hat{C}} @>>> 0,
\end{CD}
$$

\noindent the length of the skyscraper sheaf $f_{C*}\mO_{C}/\mO_{\hat{C}}$ is
3. Since a section over C can have at most two zeros, there are
 three singular points in $\hat{C}$, each of whose inverse images consists of
exactly two points, counting multiplicity.

For $C$ a smooth jumping conic with $a=1$, $E(1)_{C} \cong \Op{1}(3)\oplus \Op{1}(1).$
Consider the injective homomorphism r

$$
W\equiv H^{0}(\PS{2};E(1))\overset{rest.}{\to} H^{0}(C;E(1)_{C})
\to H^{0}(\PS{1};\Op{1}(3))\oplus H^{0}(\PS{1};\Op{1}(1))
$$

\noindent and write $r(\sigma)=(r_{1}(\sigma),r_{2}(\sigma))=(\sigma_{1},\sigma_{2}).$
$E(1)$ being globally generated implies that $r_{2}$ is
surjective. Let $\tau \in H^{0}(\PS{2};E(1))$ span $\Ker r_{2}.$
If $\sigma \in W$ is a multiple of $\tau$ (i.e. $\sigma_{2}=0$), it has three zeros on
$C.$ Otherwise the only zero of $\sigma$ on $C$ is the unique zero
of $\sigma_{2}.$ Therefore the triple point $\tau$ is the
lone singular point of $\hat{C}$ and, as before, $\hat{C}$ has
degree 4 and $l(f_{C*}\mO_{C}/\mO_{\hat{C}})=3.$

For a smooth jumping conic with $a=2$, $E(1)_{C} \cong \Op{1}(4)\oplus \Op{1}.$
Clearly $dim \{\sigma \in W \mid \Op{1}-\text{component of }
\sigma \res{C} =0\}=2$; let $\sigma$, $s$ be a basis. It follows
that $C=Z_{\sigma \wedge s}.$ Conversely, if $\sigma$, $s$ are
independent in $W$, $Z_{\sigma \wedge s}$ is a conic $C.$ It is
clear that the zero set of $\sigma$, as a point set, lies in $C.$
We want to verify that $l(Z_{\sigma})=4$ implies that $l(Z_{\sigma \text{ }\res{C}})=4$
and therefore $E(1)_{C} \cong \Op{1}(4)\oplus \Op{1}.$ For $p \in Z_{\sigma}$,
choose an affine part $\cnums ^{2}$ of $\PS{2}$ containing $p$ and
a trivialization of $E(1)$ on $\cnums ^{2}$ so that $\sigma =(\sigma_{1}, \sigma_{2})$
and $s=(s_{1}, s_{2}).$ The quadratic equation for $C$ on $\cnums ^{2}$
is $g=\sigma_{1}s_{2}-\sigma_{2}s_{1}=0.$ The
multiplicity of the zero $p$ of $\sigma$ is defined as
$i_{p}(\sigma) \equiv \text{dim }\mO_{\PS{2} ,p}/(\sigma_{1}, \sigma_{2})_{p}.$
The multiplicity of $\sigma \res{C}$ at $p$ is

\begin{align}
i_{p} &\equiv \text{dim }\mO_{C,p}/(\sigma_{1}\res{C}, \sigma_{2}\res{C})_{p}
 \notag \\
  &=\text{dim }\mO_{\PS{2} ,p}/(g)_{p} \bigg/ (\sigma_{1},
  \sigma_{2})_{p}/(g)_{p}  \notag \\
  &=\text{dim }\mO_{\PS{2} ,p}/(\sigma_{1}, \sigma_{2})_{p}
  \notag \\
  &\equiv i_{p}(\sigma). \notag
\end{align}

\noindent Remark: We have shown that a smooth jumping conic $C$ has jump
size $a=2$ iff $C=Z_{\sigma \wedge s}$ and thus $\hat{C}$ is the
line spanned by $\sigma$ and $s$ in $\mP W$ and $f_{C} \colon C \to \PS{1}$
is a 4-to-1 branched cover. There are 6 ramification points,
counting multiplicity and, as a point set, $R_{f_{C}}\subset R.$
Since $R\cdot C =6$, this suggests that $R_{f_{C}}=R\cdot C$ as
divisors on $C$ but this has not been proved.

We summarize in

\begin{Prop}  \label{Prop:(0,3)general type}
Let $E$ be a rank 2 stable bundle on $\PS{2}$ with $c_{1}=0$, $c_{2}=3$
of general type and let $f \colon \PS{2} \to \PS{2}$ be the
corresponding regular map (see (\ref{E:e(0,3)2})). Let $f_{L}$ and
$f_{C}$ be the restrictions of $f$ to a line $L$ and a smooth
conic $C$ respectively. Then \\
i) $L$ is generic iff $f_{L}$ is 1-to-1 onto a smooth conic. \\
ii) $L$ is a jumping line iff the jump size is $1$ and $f_{L}$ is
a 2-to-1 branched cover onto a line. \\
iii) $C$ is a generic conic iff $f_{C}$ is generically 1-to-1 onto
a curve of degree 4 with three singular points
each having inverse image of cardinality 2. \\
iv) $C$ is a jumping conic of jump size 1 iff $f_{C}$ is 1-to-1
onto a curve of degree 4 which is smooth except for a lone singular point $\tau$
for which $f^{-1}_{C}(\tau)$ consists of 3 points. \\
v) $C$ is a jumping conic of jump size 2 iff $C=Z_{\sigma \wedge s}$
for $\sigma$ and $s$ global sections of $E(1)$ and, in this case, $f_{C}$
is a 4-to-1 branched cover onto a line.
\end{Prop}

Now consider a bundle $E$ of non-general type meaning that there
is a line $L \subset \PS{2}$ and a $\sigma \in H^{0}(\PS{2};E(1))$
with at least 3 zeros on $L.$

\begin{Prop} \label{Prop:(0,3)non-general type}
Let $E$ be a rank 2 stable bundle on $\PS{2}$ with $c_{1}=0$, $c_{2}=3$
of non-general type. Then there is a unique line $L$ such that
every global section of $E(1)$ has exactly 3 zeros on $L.$ $E$ can
be obtained from $T\PS{2}$ by an elementary modification of the
form

\begin{equation} \label{E:e(0,3)3}
\sexs{E(1)}{T\PS{2}}{j_{L*}\mO_{L}(4)}.
\end{equation}
\end{Prop}

\begin{proof}
There is a line $L$ and $\sigma \in H^{0}(\PS{2};E(1))$ such that
$\sigma$ has 3 zeros on $L.$ All 4 of the zeros of $\sigma$ can
not be on $L$ or else the cohomology sequence of

$$
\begin{CD}
0 @>>> \Op{2}(-1) @> \sigma >> E @> \sigma \wedge >>
\mI_{Z_{\sigma}}(1) @>>> 0
\end{CD}
$$

\noindent would imply $h^{0}(\PS{2};E)=h^{0}(\PS{2};\mI_{Z_{\sigma}}(1))=1$
contradicting stability. Clearly $E(1)_{L}\cong \mO_{L}(3)\oplus \mO_{L}(-1)$
and therefore every global section of $E(1)$ has 3 zeros on $L.$
Perform the elementary modification

\begin{equation} \label{E:e(0,3)4}
\begin{CD}
0 @>>> E^{\prime} @>>> E(1) @> \phi >> j_{L*}\mO_{L}(-1) @>>> 0
\end{CD}
\end{equation}

\noindent where $\phi(s) \equiv \mO_{L}(-1)-\text{component of } s\res{L}.$
By (\ref{E:pre4}) we have $c_{1}(E^{\prime})=1$ and $c_{2}(E^{\prime})=1$ and therefore
$c_{1}(E^{\prime}(1))=3$ and $c_{2}(E^{\prime}(1))=3.$
Since $E^{\prime}$ is also stable, it follows that $E^{\prime}(1)\cong T\PS{2}$
\cite[Sec.8]{Hulek}.
By the "inverse" elementary modification \cite[page 41]{Friedman}(or by taking the
dual of (\ref{E:e(0,3)4})) one gets

$$
\begin{CD}
0 @>>> E @>>> E^{\prime} @> \psi >> j_{L*}\mO_{L}(3) @>>> 0
\end{CD}
$$

\noindent and therefore

$$
\begin{CD}
0 @>>> E(1) @>>> T\PS{2} @> \psi >> j_{L*}\mO_{L}(4) @>>> 0.
\end{CD}
$$

\noindent Since $T\PS{2}$ is globally generated, it follows that $E(1)$
is generated by global sections on $\PS{2} \setminus L.$ This
proves the uniqueness of $L.$
\end{proof}

Note that $\psi$ is given by

\begin{equation} \label{E:e(0,3)5}
T\PS{2} \to T\PS{2}\res{L} \cong \mO_{L}(2)\oplus \mO_{L}(1)
\overset{\psi_{2}\oplus \psi_{3}}{\longrightarrow} \mO_{L}(4)
\end{equation}

\noindent and so, counting parameters,
$dim \{E\}=dim \{L\}+ dim\{(\psi_{2},\psi_{3})\} -
\text{dim }Aut(\mO_{L}(4))=2+3+4-1=8$ (Recall that $\text{dim }\mM (0,3) =9.).$

We describe the jumping lines and conics in terms of the
elementary modification (\ref{E:e(0,3)3}). Clearly $L$ is a
jumping line of jump size 2. For any other line $L_{1}$, set $p=L \cap L_{1}$
and restrict (\ref{E:e(0,3)3}) to $L_{1}$ to get

$$
\begin{CD}
0 @>>> E(1)_{L_{1}} @>>> T\PS{2}_{L_{1}} @> \psi_{L_{1}} >>
\cnums p @>>> 0
\end{CD}
$$

\noindent where $\cnums p$ is $\cnums$ at $p$ and $0$ elsewhere.
Examine $\psi_{L_{1}}$ by choosing a splitting of $T\PS{2}_{L_{1}}$

\begin{equation} \label{E:e(0,3)6}
T\PS{2}_{L_{1}} \cong TL_{1}\oplus NL_{1} \cong
\mO_{L_{1}}(2)\oplus \mO_{L_{1}}(1) \to {\mO_{L_{1}}(2)}_{p}\oplus
{\mO_{L_{1}}(1)}_{p} \to \cnums.
\end{equation}

\noindent Set $K_{p}= \Ker \psi_{L_{1}}.$ From the definition of
the elementary modification (\ref{E:e(0,3)3}), one sees that $K_{p}=TL_{1p}$
iff $E(1)_{L_{1}}\cong \mO_{L_{1}}(2)\oplus \mO_{L_{1}}$ and $K_{p} \ne TL_{1p}$
iff $E(1)_{L_{1}}\cong \mO_{L_{1}}(1)\oplus \mO_{L_{1}}(1).$

Recall that the jumping line locus $J_{1} \subset \PS{2}$ is a
curve of degree 3. The lines in ${\PS{2}}^{*}$ have the form
$\Lambda_{p}\equiv \text{ lines in }\PS{2} \text{ through p}$ and
therefore, if $\Lambda_{p} \nsubseteq J_{1}$, $J_{1}\cdot \Lambda_{p}=3$,
i.e. there are three jumping
lines through each point p. If $p \in L$, the fact that the jump
size of $L$ is 2 implies that $(J_{1}\cdot \Lambda_{p})_{l}\geq 2$
for $l$ the point in ${\PS{2}}^{*}$ corresponding to $L.$ View $\psi_{p}$
via a splitting of $T\PS{2}_{L}$ as in (\ref{E:e(0,3)5}) and let $z,w \in L$
be the zeros of $\psi_{2}.$ For $p \ne w, z$, $\Ker \psi_{p} \ne
TL_{p}$ and so, by the discussion of the previous paragraph,
 there is a jumping line through $p$ transverse to $L.$
For $p=z \text{ or }w$, $(J_{1}\cdot \Lambda_{p})_{l}=3.$ This
shows that, for $z \neq w$, $l$ is a node of $J_{1}$ and that
$\Lambda_{z}$ and $\Lambda_{w}$
are the tangents to $J_{1}$ at $l.$ If $z=w$, $J_{1}$ has a cusp at $l.$

For $p \in \PS{2} \setminus L$ the jumping lines through $p$ are
determined as follows. For $a \neq b \in L$, parameterize $L$ by
$\gamma(t)=t_{0}a+t_{1}b.$ The line $\overline{p\gamma(t)}$ is
jumping iff $p \text{ mod } \cnums \gamma(t) \in T\PS{2}_{\gamma(t)}$
spans $K_{\gamma(t)}.$ Now $\xi (t)=p \text{ mod } \cnums \gamma(t)$
defines a global section of $T\PS{2}_{L}(-1)$, in fact, it gives a never vanishing
section of $NL(-1) \cong \Op{1}$ and so
$\psi_{\gamma(t)}(\xi (t))= \psi_{3}(\gamma(t)) \in H^{0}(L;\mO_{L}(3))$
is not identically zero. Its
three zeros determine the jumping lines through $p.$ We see that $l$
is the only singular point of $J_{1}$ as follows. For $L_{1}\ne L$
a jumping line and $l_{1}$ the corresponding point of ${\PS{2}}^{*}$,
let $p=L\cap L_{1}$ as above. Then $J_{1}\cdot \Lambda_{p}=2l+l_{1}$
so $(J_{1}\cdot \Lambda_{p})_{l_{1}}=1.$ If $l_{1}$ were a
singular point of $J_{1}$, $(J_{1}\cdot \Lambda)_{l_{1}}\geq 2$
for all lines through $l_{1}.$ We have proven

\begin{Prop}
Let $E$ be a rank 2 stable bundle with $c_{1}=0$, $c_{2}=3$ of
non-general type. The jumping lines form a cubic curve in ${\PS{2}}^{*}$
with exactly one singular (double) point corresponding to the line $L$ (in
Proposition~\ref{Prop:(0,3)non-general type}).
\end{Prop}

Let $C$ be a smooth conic and as usual give it a quadratic
parameterization $\PS{1}\overset{\gamma}{\to} C.$ The splittings

$$
T\PS{2}_{C} \cong \Op{1}(3)\oplus \Op{1}(3)\cong U\otimes \Op{1}(3)
$$

\noindent for $\text{dim }U=2$ follow from examination of the global
sections of $T\PS{2}.$ We determine if $C$ is a jumping conic by
noting the effect of the elementary modification (\ref{E:e(0,3)3})
on this splitting. Let $p$ and $q$ be the intersection points of $C$
and $L$ and assume they are distinct. Restricting
(\ref{E:e(0,3)3}) to $C$, we get

$$
\begin{CD}
0 @>>> E(1)_{C} @>>> T\PS{2}_{C} @> \tilde{\psi_{p}}\oplus \tilde{\psi_{q}}
>> \cnums p \oplus \cnums q @>>> 0
\end{CD}
$$

\noindent where $\tilde{\psi_{p}}$ is obtained from the homomorphism $\psi$
of (\ref{E:e(0,3)5}) evaluated at $p$ and similarly for $\tilde{\psi_{q}}.$
Therefore $K_{p}\equiv \Ker \tilde{\psi_{p}}$ and $K_{q}\equiv \Ker \tilde{\psi_{q}}$
can be considered as 1-dimensional subspaces of $U.$ If $K_{p}=K_{q}$, $E(1)_{C}$
contains a copy of $\Op{1}(3)$ and so $E(1)_{C}\cong \Op{1}(3)\oplus \Op{1}(1)$.
If $K_{p} \neq K_{q}$, $E(1)_{C}$
contains a copy of $\Op{1}(2)\oplus \Op{1}(2)$ and so has this
splitting.

Now suppose $C$ is tangent to $L$ at $p.$ By a change in the homogeneous
coordinates of $\PS{2}$ we can assume that $p=e_{0}$ ($p=(0,0)$ in the
corresponding affine coordinates $x_{1}$, $x_{2}$) and $L$ is
defined by $x_{2}=0.$ By changing coordinates on a smaller Zariski
open neighborhood of $p$ we can also assume that $C$ is defined by
$x_{2}+x_{1}^{2}=0$ and therefore $x_{1}$ is a local coordinate of $C$
near $p.$ Restricting (\ref{E:e(0,3)3}) to $C$ we have an
elementary modification at the divisor $2e_{0}$

$$
\begin{CD}
0 @>>> E(1)_{C} @>>> T\PS{2}_{C} @> \Psi
>> \cnums[x_{1}]/(x_{1}^{2}) @>>> 0
\end{CD}
$$

\noindent where
$\Psi \colon T\PS{2}_{C}\cong U\otimes \Op{1}(3)\to \cnums[x_{1}]/(x_{1}^{2})$
is determined by the one-jet $j_{1}(\psi)$ of $\psi$ at p (see
(\ref{E:e(0,3)5}))
as follows. With respect to a basis for $U$ and the dual basis for
$U^{*}$ and for $f$ a local section of $T\PS{2}_{C}$,

$$
f(x_{1})=\left(\begin{matrix}
f_{1}(x_{1}) \\
f_{2}(x_{1})
\end{matrix} \right),
j_{1}(\psi)=\left( \begin{matrix}
a_{10}+a_{11}x_{1} \\
a_{20}+a_{21}x_{1}
\end{matrix} \right)
$$

\noindent and

\begin{equation} \label{E:e(0,3)7}
\Psi(f)=j_{1}((a_{10}+a_{11}x_{1})f_{1}(x_{1})+(a_{20}+a_{21}x_{1})f_{2}(x_{1})).
\end{equation}

We show that $C$ is a jumping conic if and only if the 1-jets $a_{10}+a_{11}x_{1}$
and $a_{20}+a_{21}x_{1}$ are scalar multiples iff $a_{10}a_{21}-a_{11}a_{20}=0$
iff $\text{ker }\psi(0)=\text{ker }\psi^{\prime}(0)$ as subspaces
of $U.$ To verify this, change basis in $U$ so that $a_{10}=0$ and
$a_{20}=1$ and therefore the above condition is $a_{11}=0.$
Let $t_{0}$, $t_{1}$ be homogeneous coordinates of a
quadratic parameterization $\PS{1} \to C$ with $[1,0]$ corresponding to p,
$t \equiv t_{1}/t_{0}$
the affine coordinate on $\mU_{0} \equiv \PS{1}\setminus \infty$,
and $s \equiv 1/t$ the coordinate on $\mU_{1} \equiv \PS{1}\setminus 0.$
Then $f \in T\PS{2}$ is in $E(1)$ iff $\Psi (f)=0$ iff

\begin{align} \label{E:e(0,3)8}
f_{2}(0)&=0     \\
f_{2}^{\prime}(0)+a_{11} f_{1}(0)&=0.   \notag
\end{align}

\noindent Assume $\text{ker }\psi(0)=\text{ker }\psi^{\prime}(0)$, that is, $a_{11}=0.$
Then

$$
f(t)=\left( \begin{matrix}
f_{1}(t)  \\
t^{2}h_{2}(t)
\end{matrix} \right)
$$

\noindent and this implies $E(1)_{C} \cong \Op{1}(3)\oplus \Op{1}(1).$
Assuming instead that $a_{11} \ne 0$,

\begin{equation} \label{E:e(0,3)9}
f(t)=\left( \begin{matrix}
-a_{11}^{-1}f_{21}+f_{11}t+t^{2}h_{1}(t)  \\
f_{21}t+t^{2}h_{2}(t)  \notag
\end{matrix} \right)
=\left( \begin{matrix}
t & -a_{11}^{-1} \\
0 & t
\end{matrix} \right) \cdot \left( \begin{matrix}
k_{1}(t) \\
k_{2}(t)
\end{matrix} \right)
\end{equation}

\noindent where $k_{1}(t)=f_{11}+th_{1}(t)+a_{11}^{-1}h_{2}(t)$
and $k_{2}(t)=f_{21}+th_{2}(t).$ The transition function for $E(1)$
(from the $\mU_{0}$ to the $\mU_{1}$ trivialization) is therefore

$$
\left( \begin{matrix}
t^{-3}&0  \\
0& t^{-3}
\end{matrix} \right) \cdot \left( \begin{matrix}
t&-a_{11}^{-1} \\
0&t
\end{matrix} \right) = \left( \begin{matrix}
t^{-2}& -a_{11}^{-1} t^{-3} \\
0& t^{-2}
\end{matrix} \right) = \left( \begin{matrix}
1&-a_{11}^{-1}s \\
0&1
\end{matrix}  \right) \cdot \left( \begin{matrix}
t^{-2}&0 \\
0& t^{-2}
\end{matrix} \right).
$$

\noindent Changing the trivialization of $E(1)$ over $\mU_{1}$ via
the inverse of the matrix function of $s$ above, we see that
$E(1)_{C} \cong \Op{1}(2)\oplus \Op{1}(2).$

\section{Rank 2 Stable Bundles on $\PS{2}$ with $c_{1}=-1$ and
$c_{2}=2$} \label{S:(-1,2)}

The analysis of jumping lines and conics for a rank 2 stable
bundle $E$ on $\PS{2}$ with $c_{1}=-1$ and $c_{2}=2$ is very similar
to the study of the $c_{1}=0$, $c_{2}=3$ non-general type case. We
summarize the results. $E$ has a unique jumping line $L$
characterized by: $L=Z_{\sigma \wedge s}$ for $\sigma$, $s$ a
basis for $H^{0}(\PS{2};E(1)).$ $E$ can be expressed as an
elementary modification along $L$ of the trivial bundle

\begin{equation} \label{E:e(-1,2)1}
\begin{CD}
0 @>>> E @>>> U \otimes \Op{2} @> \psi >>
j_{L*}\mO_{L}(2) @>>> 0
\end{CD}
\end{equation}

\noindent for $U$ a vector space of dimension 2. A smooth conic $C$
which intersects $L$ at distinct points $p$ and $q$ is a jumping
conic iff $\text{ker }\psi_{p}=\text{ker }\psi_{q}$ as subspaces
of $U.$ If $C$ and $L$ have a double intersection at $p$, $C$ is a
jumping conic iff $\text{ker }\psi_{p}=\text{ker }\psi^{\prime}_{p}.$
Furthermore, if homogeneous coordinates are chosen so that $L$ is
defined by $x_{2}=0$ and if $\psi$ is given by
$\psi_{k}=\psi_{k00}x_{0}^{2}+\psi_{k01}x_{0}x_{1}+\psi_{k11}x_{1}^{2}$
for $k=1,2$, then the hyperplane $J_{2}\subset \PS{5}$ is
given by

$$
a_{00}\xi_{00}+a_{01}\xi_{01}+a_{11}\xi_{11}=0
$$

\noindent where $C$ is defined by $\sum_{0 \leq i \leq j \leq 2}\xi_{ij}x_{i}x_{j}=0$
and $a_{00}=\psi_{101}\psi_{211}-\psi_{111}\psi_{201}$,
$a_{01}=\psi_{111}\psi_{200}-\psi_{100}\psi_{211}$, and
$a_{11}=\psi_{100}\psi_{201}-\psi_{101}\psi_{200}.$

\providecommand{\bysame}{\leavevmode\hbox to3em{\hrulefill}\thinspace}
\providecommand{\MR}{\relax\ifhmode\unskip\space\fi MR }
\providecommand{\MRhref}[2]{%
  \href{http://www.ams.org/mathscinet-getitem?mr=#1}{#2}
}
\providecommand{\href}[2]{#2}


\begin{thebibliography}{10}

\bibitem{Barth-P2}
W.~Barth, \emph{Moduli of vector bundles on the projective plane}, Invent.
  Math. \textbf{42} (1977), 63--91. \MR{57 \#324}

\bibitem{Barth-sb}
\bysame, \emph{Some properties of stable rank-$2$ vector bundles on ${\bf
  {P}}\sb{n}$}, Math. Ann. \textbf{226} (1977), no.~2, 125--150. \MR{55 \#2905}

\bibitem{Hartshorne-rt}
Lawrence Ein, Robin Hartshorne, and Hans Vogelaar, \emph{Restriction theorems
  for stable rank $3$ vector bundles on ${\bf {P}}\sp{n}$}, Math. Ann.
  \textbf{259} (1982), no.~4, 541--569. \MR{84b:14006}

\bibitem{Friedman}
Robert Friedman, \emph{Algebraic surfaces and holomorphic vector bundles},
  Springer-Verlag, New York, 1998. \MR{99c:14056}

\bibitem{GH}
Phillip Griffiths and Joseph Harris, \emph{Principles of algebraic geometry},
  Wiley Classics Library, John Wiley \& Sons Inc., New York, 1994, Reprint of
  the 1978 original. \MR{95d:14001}

\bibitem{splitting}
A.~Grothendieck, \emph{Sur la classification des fibr\'es holomorphes sur la
  sph\`ere de {R}iemann}, Amer. J. Math. \textbf{79} (1957), 121--138.
  \MR{19,315b}

\bibitem{Hartshorne}
Robin Hartshorne, \emph{Algebraic geometry}, Springer-Verlag, New York, 1977,
  Graduate Texts in Mathematics, No. 52. \MR{57 \#3116}

\bibitem{Hartshorne-sb}
\bysame, \emph{Stable vector bundles of rank $2$ on ${\bf {P}}\sp{3}$}, Math.
  Ann. \textbf{238} (1978), no.~3, 229--280. \MR{80c:14011}

\bibitem{Hartshorne-srs}
\bysame, \emph{Stable reflexive sheaves}, Math. Ann. \textbf{254} (1980),
  no.~2, 121--176. \MR{82b:14011}

\bibitem{Hulek}
Klaus Hulek, \emph{Stable rank-{$2$} vector bundles on {${\bf P}\sb{2}$} with
  {$c\sb{1}$}\ odd}, Math. Ann. \textbf{242} (1979), no.~3, 241--266.
  \MR{80m:14011}

\bibitem{Hulek-hm1}
\bysame, \emph{Geometry of the {H}orrocks-{M}umford bundle}, Algebraic
  geometry, Bowdoin, 1985 (Brunswick, Maine, 1985), Proc. Sympos. Pure Math.,
  vol.~46, Amer. Math. Soc., Providence, RI, 1987, pp.~69--85. \MR{89a:14017}

\bibitem{Hurtubise-jumpinglines}
Jacques Hurtubise, \emph{Instantons and jumping lines}, Comm. Math. Phys.
  \textbf{105} (1986), no.~1, 107--122. \MR{87g:14009}

\bibitem{Manaresi}
Mirella Manaresi, \emph{On the jumping conics of a semistable rank two vector
  bundle on {${\bf P}\sp 2$}}, Manuscripta Math. \textbf{69} (1990), no.~2,
  133--151. \MR{92b:14023}

\bibitem{Maruyama1}
M.~Maruyama, \emph{Boundedness of semistable sheaves of small ranks}, Nagoya
  Math. J. \textbf{78} (1980), 65--94.

\bibitem{Maruyama2}
Masaki Maruyama, \emph{Moduli of stable sheaves. {I}{I}}, J. Math. Kyoto Univ.
  \textbf{18} (1978), no.~3, 557--614. \MR{82h:14011}

\bibitem{OSS}
Christian Okonek, Michael Schneider, and Heinz Spindler, \emph{Vector bundles
  on complex projective spaces}, Birkh\"auser Boston, Mass., 1980.
  \MR{81b:14001}

\bibitem{Ran-jumping}
Ziv Ran, \emph{The degree of the divisor of jumping rational curves}, Q. J.
  Math. \textbf{52} (2001), no.~3, 367--383. \MR{2002j:14009}

\bibitem{Schwarzenberger1}
R.~L.~E. Schwarzenberger, \emph{Vector bundles on the projective plane}, Proc.
  London Math. Soc. (3) \textbf{11} (1961), 623--640. \MR{25 \#1161}

\end{thebibliography}
\end{document}